\documentclass[12pt]{article}
\usepackage{latexsym,amsfonts,amssymb,mathrsfs}
\usepackage[dvips]{color}
\usepackage{colordvi,multicol}

\def \cal{\mathcal}

\newtheorem{thm}{Theorem}[section]

\newtheorem{lem}[thm]{Lemma}

\begin{document}
\title{\bf The 6-element case of S-Frankl conjecture (I)  }

\author {
  Ze-Chun Hu and Shi-Lun Li\thanks{Corresponding author: College of Mathematics, Sichuan University, 610064, China\vskip 0cm E-mail addresses: lishilun@scu.edu.cn (S.-L. Li)}\\ \\
   College of Mathematics, Sichuan University, China}

\maketitle
\date{}


\vskip 0.5cm \noindent{\bf Abstract}\quad  The union-closed sets conjecture (Frankl's conjecture)  says that  for any finite union-closed family of finite sets, other than the family consisting only of
the empty set, there exists an element that belongs to at least half of the sets in the family. In \cite{CH17},  a stronger version of Frankl's conjecture
(S-Frankl conjecture for short) was introduced and  a partial proof was given. In particular, it was proved in \cite{CH17} that S-Frankl conjecture holds when $n\leq 5$, where  $n$ is the number of all the elements in the family of sets. Now, we want to prove that it holds when $n=6$.  Since the paper is very long, we split it into two parts. This is the first part.

\smallskip

\smallskip



\section{Introduction}

A family $\mathcal{F}$ of sets is union-closed if  $A,B\in \mathcal{F}$ implies
 $A\cup  B\in \mathcal{F}$. For simplicity, denote $n=|\cup_{A\in \cal{F}}A|$ and $m=|\cal{F}|$.

In 1979, Peter Frankl (cf. \cite{Ri85, St86}) conjectured that  for any finite union-closed family of finite sets, other than the family consisting only of the empty set, there exists an element that belongs to at least half of the sets in the family.

If a union-closed family $\cal{F}$ contains a set with one element or two elements, then Frankl's conjecture holds for $\cal{F}$ (\cite{SR89} ).
The result was extended by Poonen (\cite{Po92}). In addition,  the author in \cite{Po92}  proved that Frankl's conjecture holds if $n\leq 7$ or $m\leq 28$, and presented an equivalent lattice formulation of Frankl's conjecture. Bo\v{s}njak and Markovi\'{c}  proved  that Frankl's conjecture holds if $n\leq 11$.   Zivkovi\'{c} and Vu\v{c}kovi\'{c}  gave a computer
assisted proof   in ``The 12 element case of Frankl's conjecture, Preprint (2012)"that Frankl's conjecture holds if  $n\leq 12$, which together with Faro's result (\cite{Lo94-b}) (see also Roberts and Simposon \cite{RS10}) implies that Frankl's conjecture is true if  $m\leq 50$.  For more progress on Frankl's conjecture, we refer to \cite{BS15}, \cite{GY98}, \cite{JV98}, \cite{Lo94-a}, \cite{Ma07}, \cite{Mo06},   \cite{SR89}, \cite{Va02}, \cite{Va03}, \cite{Va04}.

Let $M_n=\{1,2,\ldots,n\}$ and $\cal{F}\subset 2^{M_n}=\{A: A\subset M_n\}$ with $\cup_{A\in\cal{F}}A=M_n$.  Suppose that $\cal{F}$ is union-closed. Without loss of generality, we assume that $\emptyset \in \cal{F}$.
For any $k=1,2,\ldots,n$, define
$
\cal{M}_k=\{A\in 2^{M_n}: |A|=k\},
$
and
$$
T(\cal{F})=\inf\{1\leq k\leq n: \cal{F}\cap \cal{M}_k\neq \emptyset\}.
$$
Then $1\leq T(\cal{F})\leq n$.  By virtue of $T(\cal{F})$, the authors in \cite{CH17}  introduced the following stronger version of Frankl's conjecture (S-Frankl conjecure for short).

{\bf S-Frankl conjecture:} If  $n\geq 2$ and $T(\cal{F})=k\in \{2,\ldots,n\}$, then there exist at least $k$ elements in $M_n$ which belong to at least half of the sets in $\cal{F}$.

%
%
%

We  need the following lemma, which has been used in some proofs in \cite{CH17}.

\begin{lem}\label{lem-1.2}
Suppose that $M$ is a finite set with $|M|\geq 2$ and $\cal{G}\subset \{A\subset  M: |A|=|M|-1\}$. If $|\cal{G}|\geq 2$, then all the elements in $M$ belong to at least $|\cal{G}|-1$ set(s) in $\cal{G}$.
\end{lem}
{\bf Proof.} Without loss of generality, we assume that $M=M_n=\{1,2,\ldots,n\}$ with $n\geq 2$. Then $\cal{G}$ is a subset of $\cal{M}_{n-1}=\{A\subset M_n: |A|=n-1\}$. Notice that for any $i\in \{1,2,\ldots,n\}$, it belongs to all the sets in $\cal{M}_{n-1}$ except the set $M_n\backslash \{i\}$. Hence all the elements in $M$ belong to at least $|\cal{G}|-1$ set(s) in $\cal{G}$.\hfill\fbox

When we consider S-Frankl conjecture for the case that $n=6$, by Section 2 of \cite{CH17}, we know that if $T(\cal{F})\in \{4, 5,6\}$, then there exist at least $T(\cal{F})$ elements in $M_6$ which belong to at least half of the sets in $\cal{F}$.
Thus we need only to consider the two cases $T(\cal{F})=3$ and  $T(\cal{F})=2$.  In next section, we will prove that S-Frankl conjecture holds when $n=6$ and $T(\cal{F})=3$. The proof for the case that $n=6$ and $T(\cal{F})=2$ will be given in a sister paper.

\section{S-Frankl conjecture for $n=6$ and $T(\cal{F})=3$}


 Let $M_6=\{1,2,\ldots,6\}$ and $\cal{F}\subset 2^{M_6}=\{A: A\subset M_6\}$ with $\cup_{A\in\cal{F}}A=M_6$.  Suppose that $\cal{F}$ is union-closed and  $\emptyset \in \cal{F}$.  For $k=1,2,\ldots,6$, define $
\cal{M}_k=\{A\in 2^{M_6}: |A|=k\}, n_k=|\cal{F}\cap \cal{M}_k|,
$ and
$$
T(\cal{F})=\inf\{1\leq k\leq 6: n_k>0\}.
$$
Then $1\leq T(\cal{F})\leq 6$.

In the following, we assume that $T(\cal{F})=3$, and   will prove that there exist at least 3 elements in $M_6$ which belong to at least half of the sets in $\cal{F}$.
We have  4 cases:   $\cal{F}=\{\emptyset, M_6\}\cup \cal{G}_3$,  $\cal{F}=\{\emptyset, M_6\}\cup \cal{G}_3\cup \cal{G}_5$,  $\cal{F}=\{\emptyset, M_6\}\cup \cal{G}_3\cup \cal{G}_4$ and  $\cal{F}=\{\emptyset, M_6\}\cup \cal{G}_3\cup \cal{G}_4\cup \cal{G}_5$, where $\cal{G}_i$ is a nonempty subset of $\cal{M}_i$ for $i=3,4,5$.

\subsection{$\cal{F}=\{\emptyset, M_6\}\cup \cal{G}_3$}

We have two subcases: $n_3=1$ and $n_3\geq 2$. Throughout  the rest of this paper, we omit the sentences of this type.  Denote $\cal{G}_3=\{G_1,\ldots,G_{n_3}\}$.

\begin{itemize}
\item[(1)] $n_3=1$. Now $\cal{G}_3=\{G_1\}$. Then all the 3 elements in $G_1$ belong to two sets among the three sets in $\cal{F}$.

\item[(2)]  $n_3\geq 2$. For any $i,j=1,\ldots,n_3,i\neq j$, we must have $G_i\cup G_j=M_6$, which implies that $G_i\cap G_j=\emptyset$. Hence $n_3=2$. Now all the 6
elements in $M_6$ belong to two sets among the four sets in $\cal{F}$.
\end{itemize}

\subsection{$\cal{F}=\{\emptyset, M_6\}\cup \cal{G}_3\cup \cal{G}_5$}

Denote $\cal{G}_3=\{G_1,\ldots,G_{n_3}\}$ and $\cal{G}_5=\{H_1,\ldots,H_{n_5}\}$.

\begin{itemize}
\item[(1)] $n_5=1$.  Now $\cal{G}_5=\{H_1\}$.  Without loss of generality, we assume that $H_1=\{1,2,3,4,5\}$.
\begin{itemize}
\item[(1.1)] $n_3=1$. Now $\cal{G}_3=\{G_1\}$. Notice that  all the elements in $G_1\cup \{1,2,3,4,5\}$ belong to at least one of the two sets $G_1$ and $\{1,2,3,4,5\}$. Then we know that  all the elements in $G_1\cup \{1,2,3,4,5\}$ belong to at least half of the sets in $\cal{F}$.

\item[(1.2)] $n_3\geq 2$. For any $i,j=1,\ldots,n_3, i\neq j$, we have $G_i\cup G_j=M_6$ or $G_i\cup G_j=\{1,2,3,4,5\}$.

\begin{itemize}
\item[(1.2.1)] $n_3$ is an even number and there is a permutation $(i_1,\ldots,i_{n_3})$ of $(1,\ldots,n_3)$ such that $G_{i_1}\cup G_{i_2}=\cdots=G_{i_{n_3-1}}\cup G_{i_{n_3}}=M_6$. Then all the 6 elements in
$M_6$ belong to half of the sets in $\cal{G}_3$ and thus all the 5 elements in $\{1,2,3,4,5\}$ belong to at least half of the sets in $\cal{F}$.

\item[(1.2.2)] $n_3$ is an odd number and there is a permutation $(i_1,\ldots,i_{n_3})$ of $(1,\ldots,n_3)$ such that $G_{i_1}\cup G_{i_2}=\cdots=G_{i_{n_3-2}}\cup G_{i_{n_3-1}}=M_6$. Then all the 6 elements in
$M_6$ belong to half of the sets in $\{G_{i_1},\ldots,G_{i_{n_3-1}}\}$ and thus all the  elements in $G_{i_{n_3}}\cup \{1,2,3,4,5\}$ belong to at least half of the sets in $\cal{F}$.

\item[(1.2.3)] We can decompose $\cal{G}_3$ into two disjoint parts $\{G_{i_1},\ldots,G_{i_{2k}}\}$ (hereafter this part may be an empty set) and $\{G_{i_{2k+1}},\ldots,G_{i_{n_3}}\}$, where
$\{i_1,\ldots,i_{n_3}\}=\{1,\ldots,n_3\}$, $n_3-2k\geq 2$, and

\quad (i) $G_{i_1}\cup G_{i_2}=\cdots=G_{i_{2k-1}}\cup G_{i_{2k}}=M_6$;

\quad (ii) for any two different indexes $\{i,j\}$ from $\{i_{2k+1},\ldots,i_{n_3}\}, G_i\cup G_j=\{1,2,3,4,5\}$.

Then all the 6 elements in $M_6$ belong to half of the sets in $\{G_{i_1},\ldots,G_{i_{2k}}\}$. Without loss of generality, we assume that $G_{i_{2k+1}}=\{1,2,3\}$. By (ii), we know that for any $j=i_{2k+2},\ldots,i_{n_3}$,
$$
G_j\in \{\{1,4,5\}, \{2, 4,5\}, \{3, 4,5\}\}.
$$
Since $|\{1,4,5\}\cup \{2,4,5\}|=|\{1,4,5\}\cup \{3,4,5\}|=|\{2,4,5\}\cup \{3,4,5\}|=4$, by (ii) again, we know that in this case $n_3-2k=2$. Then we know that all the 5 elements in
$\{1,2,3,4,5\}$ belong to at least half of the sets in $\cal{F}$.
\end{itemize}
\end{itemize}

\item[(2)] $n_5\geq 2$.  By Lemma \ref{lem-1.2}, we know that all the 6 elements in $M_6$ belong to at least $n_5-1$ set(s) in $\cal{G}_5$ and thus belong to at least half of the sets
in $\cal{G}_5$. Hence it is enough to show that there exist at least 3 elements in $M_6$ which belong to at least half of the sets in $\cal{G}_3$.

\begin{itemize}
\item[(2.1)]  $n_3=1$. Now $\cal{G}_3=\{G_1\}$. Then all the 3 elements in $G_1$ satisfy the condition.

\item[(2.2)]  $n_3\geq 2$. For any $i,j=1,\ldots,n_3,i\neq j$, we have $G_i\cup G_j=M_6$ or $|G_i\cup G_j|=5$.
\begin{itemize}
\item[(2.2.1)] $n_3$ is an even number and there is a permutation $(i_1,\ldots,i_{n_3})$ of $(1,\ldots,n_3)$ such that $G_{i_1}\cup G_{i_2}=\cdots=G_{i_{n_3-1}}\cup G_{i_{n_3}}=M_6$. Then all the 6 elements in
$M_6$ belong to half of the sets in $\cal{G}_3$.

\item[(2.2.2)] $n_3$ is an odd number and there is a permutation $(i_1,\ldots,i_{n_3})$ of $(1,\ldots,n_3)$ such that $G_{i_1}\cup G_{i_2}=\cdots=G_{i_{n_3-2}}\cup G_{i_{n_3-1}}=M_6$. Then all the 6 elements in
$M_6$ belong to half of the sets in $\{G_{i_1},\ldots,G_{i_{n_3-1}}\}$ and thus all the 3  elements in $G_{i_{n_3}}$ belong to at least half of the sets in $\cal{G}_3$.

\item[(2.2.3)] We can decompose $\cal{G}_3$ into two disjoint parts $\{G_{i_1},\ldots,G_{i_{2k}}\}$ and $\{G_{i_{2k+1}},\ldots,G_{i_{n_3}}\}$, where
$\{i_1,\ldots,i_{n_3}\}=\{1,\ldots,n_3\}$, $n_3-2k\geq 2$, and

\quad (i) $G_{i_1}\cup G_{i_2}=\cdots=G_{i_{2k-1}}\cup G_{i_{2k}}=M_6$;

\quad (ii) for any two different indexes $\{i,j\}$ from $\{i_{2k+1},\ldots,i_{n_3}\}, |G_i\cup G_j|=5$.

Then all the 6 elements in $M_6$ belong to half of the sets in $\{G_{i_1},\ldots,G_{i_{2k}}\}$. Without loss of generality, we assume that $G_{i_{2k+1}}=\{1,2,3\}$. By (ii), we know that for any $j=i_{2k+2},\ldots,i_{n_3}$,
$$
G_j\in \{\{1,4,5\}, \{2, 4,5\}, \{3, 4,5\},\{1,4,6\}, \{2, 4,6\}, \{3, 4,6\},\{1,5,6\}, \{2, 5,6\}, \{3, 5,6\}\}.
$$
Since $|\{1,4,5\}\cup \{2,4,5\}|=|\{1,4,5\}\cup \{3,4,5\}|=|\{2,4,5\}\cup \{3,4,5\}|=4$, by (ii) again, we know that
$$
|\{G_{i_{2k+1}},\ldots,G_{i_{n_3}}\}\cap \{\{1,4,5\}, \{2, 4,5\}, \{3, 4,5\}\}|\leq 1.
$$
Similarly, we have that
\begin{eqnarray*}
&&|\{G_{i_{2k+1}},\ldots,G_{i_{n_3}}\}\cap \{\{1,4,6\}, \{2, 4,6\}, \{3, 4,6\}\}|\leq 1,\\
&&|\{G_{i_{2k+1}},\ldots,G_{i_{n_3}}\}\cap \{\{1,5,6\}, \{2, 5,6\}, \{3, 5,6\}\}|\leq 1.
\end{eqnarray*}
Hence we need only to consider the following 3 cases:
\begin{itemize}
\item[(2.2.3.1)] $n_3-2k=2$. Take $
|\{G_{i_{2k+1}},\ldots,G_{i_{n_3}}\}\cap \{\{1,4,5\}, \{2, 4,5\}, \{3, 4,5\}\}|=1$ for example.
Without loss of generality, we assume that $\{G_{i_{2k+1}},\ldots,G_{i_{n_3}}\}=\{\{1,2,3\}, \{1,4,5\}\}$. Now all the 5 elements in $\{1,2,3,4,5\}$ belong to at least half of the sets in $\cal{G}_3$.

\item[(2.2.3.2)] $n_3-2k=3$. Take $
|\{G_{i_{2k+1}},\ldots,G_{i_{n_3}}\}\cap \{\{1,4,5\}, \{2, 4,5\}, \{3, 4,5\}\}|=1=|\{G_{i_{2k+1}},\ldots,G_{i_{n_3}}\}\cap \{\{1,4,6\}, \{2, 4,6\}, \{3, 4,6\}\}|$ for example.  Without loss of generality, we assume that $\{G_{i_{2k+1}},\ldots,G_{i_{n_3}}\}=\{\{1,2,3\}, \{1,4,5\},\{2,4,6\}\}$.  Now all the 3 elements in $\{1,2,4\}$ belong to at least half of the sets in $\cal{G}_3$.

\item[(2.2.3.3)] $n_3-2k=4$. Without loss of generality, we  assume that $\{G_{i_{2k+1}},\ldots,G_{i_{n_3}}\}=\{\{1,2,3\}, \{1,4,5\},\{2,4,6\}, \{3,5,6\}\}$.  Now all the 6 elements in $M_6$ belong to at least half of the sets in $\cal{G}_3$.
\end{itemize}
\end{itemize}
\end{itemize}
\end{itemize}

\subsection{$\cal{F}=\{\emptyset, M_6\}\cup \cal{G}_3\cup \cal{G}_4$}

Denote $\cal{G}_3=\{G_1,\ldots,G_{n_3}\}$ and $\cal{G}_4=\{H_1,\ldots,H_{n_4}\}$.

\begin{itemize}
\item[(1)] $n_4=1$.  Now $\cal{G}_4=\{H_1\}$. Without loss of generality, we assume that $H_1=\{1,2,3,4\}$.

\begin{itemize}
\item[(1.1)] $n_3=1$. Now $\cal{G}_3=\{G_1\}$.  Then all the elements in $G_1\cup \{1,2,3,4\}$ belong to at least half of the sets in $\cal{F}$.

\item[(1.2)] $n_3\geq 2$. For any $i,j=1,\ldots,n_3,i\neq j$, we have $G_i\cup G_j=M_6$ or $G_i\cup G_j=\{1,2,3,4\}$.
\begin{itemize}
\item[(1.2.1)] $n_3$ is an even number and there is a permutation $(i_1,\ldots,i_{n_3})$ of $(1,\ldots,n_3)$ such that $G_{i_1}\cup G_{i_2}=\cdots=G_{i_{n_3-1}}\cup G_{i_{n_3}}=M_6$. Then all the 6 elements in
$M_6$ belong to half of the sets in $\cal{G}_3$ and thus all the 4 elements in $\{1,2,3,4\}$ belong to at least half of the sets in $\cal{F}$.

\item[(1.2.2)] $n_3$ is an odd number and there is a permutation $(i_1,\ldots,i_{n_3})$ of $(1,\ldots,n_3)$ such that $G_{i_1}\cup G_{i_2}=\cdots=G_{i_{n_3-2}}\cup G_{i_{n_3-1}}=M_6$. Then all the 6 elements in
$M_6$ belong to half of the sets in $\{G_{i_1},\ldots,G_{i_{n_3-1}}\}$ and thus all the  elements in $G_{i_{n_3}}\cup \{1,2,3,4\}$ belong to at least half of the sets in $\cal{F}$.

\item[(1.2.3)] We can decompose $\cal{G}_3$ into two disjoint parts $\{G_{i_1},\ldots,G_{i_{2k}}\}$ and $\{G_{i_{2k+1}},\ldots,G_{i_{n_3}}\}$, where
$\{i_1,\ldots,i_{n_3}\}=\{1,\ldots,n_3\}$, $n_3-2k\geq 2$, and

\quad (i) $G_{i_1}\cup G_{i_2}=\cdots=G_{i_{2k-1}}\cup G_{i_{2k}}=M_6$;

\quad (ii) for any two different indexes $\{i,j\}$ from $\{i_{2k+1},\ldots,i_{n_3}\}, G_i\cup G_j=\{1,2,3,4\}$.

Then all the 6 elements in $M_6$ belong to half of the sets in $\{G_{i_1},\ldots,G_{i_{2k}}\}$. By Lemma \ref{lem-1.2}, we know that all the 4 elements in $\{1,2,3,4\}$ belong to at least $n_3-2k-1$ set(s) in $\{G_{i_{2k+1}},\ldots,G_{i_{n_3}}\}$ and thus belong to at least half of the sets in $\{G_{i_{2k+1}},\ldots,G_{i_{n_3}}\}$. Hence in this case, all the 4 elements in $\{1,2,3,4\}$ belong to at least half of the sets in $\cal{F}$.
\end{itemize}
\end{itemize}

\item[(2)] $n_4\geq 2$.  For any $i,j=1,\ldots,n_4,i\neq j$, we must have $H_i\cup H_j=M_6$. We claim that all the 6 elements in $M_6$ belong to at least half of the sets in $\cal{G}_4$.

   \quad  In fact, if $n_4=2k$ is an even number, then $H_1\cup H_2=\ldots=H_{2k-1}\cup H_{2k}=M_6$ and thus all the 6 elements in $M_6$ belong to at least half of the sets in $\cal{G}_4$. If $n_4=2k+1$ is an odd number, then by
    \begin{eqnarray*}
    &&H_1\cup H_2=\cdots=H_{2k-1}\cup H_{2k}=M_6,
    \end{eqnarray*}
    we know that all the 4 elements in $H_{2k+1}$ belong to at least half of the sets in $\cal{G}_4$;
     by
    \begin{eqnarray*}
    H_2\cup H_3=\cdots=H_{2k}\cup H_{2k+1}=M_6,
    \end{eqnarray*}
    we know that all the 4 elements in $G_1$ belong to at least half of the sets in $\cal{G}_4$. By $H_1\cup H_{2k+1}=M_6$, we know that all the 6 elements in $M_6$ belong to at least half of the sets in $\cal{G}_4$.

\quad Hence it is enough to show that there exist 3 elements in $M_6$ which belong to at least half of the sets in $\cal{G}_3$.

\begin{itemize}
\item[(2.1)]  $n_3=1$. Now $\cal{G}_3=\{G_1\}$. Then all the 3 elements in $G_1$ satisfy the condition.

\item[(2.2)]  $n_3\geq 2$. For any $i,j=1,\ldots,n_3,i\neq j$, we have $G_i\cup G_j=M_6$ or $|G_i\cup G_j|=4$.

\begin{itemize}
\item[(2.2.1)] $n_3$ is an even number and there is a permutation $(i_1,\ldots,i_{n_3})$ of $(1,\ldots,n_3)$ such that $G_{i_1}\cup G_{i_2}=\cdots=G_{i_{n_3-1}}\cup G_{i_{n_3}}=M_6$. Then all the 6 elements in
$M_6$ belong to half of the sets in $\cal{G}_3$.

\item[(2.2.2)] $n_3$ is an odd number and there is a permutation $(i_1,\ldots,i_{n_3})$ of $(1,\ldots,n_3)$ such that $G_{i_1}\cup G_{i_2}=\cdots=G_{i_{n_3-2}}\cup G_{i_{n_3-1}}=M_6$. Then all the 6 elements in
$M_6$ belong to half of the sets in $\{G_{i_1},\ldots,G_{i_{n_3-1}}\}$ and thus all the 3  elements in $G_{i_{n_3}}$ belong to at least half of the sets in $\cal{G}_3$.

\item[(2.2.3)] We can decompose $\cal{G}_3$ into two disjoint parts $\{G_{i_1},\ldots,G_{i_{2k}}\}$ and $\{G_{i_{2k+1}},\ldots,G_{i_{n_3}}\}$, where
$\{i_1,\ldots,i_{n_3}\}=\{1,\ldots,n_3\}$, $n_3-2k\geq 2$, and

\quad (i) $G_{i_1}\cup G_{i_2}=\cdots=G_{i_{2k-1}}\cup G_{i_{2k}}=M_6$;

\quad (ii) for any two different indexes $\{i,j\}$ from $\{i_{2k+1},\ldots,i_{n_3}\}, |G_i\cup G_j|=4$.

Then all the 6 elements in $M_6$ belong to  half of the sets in $\{G_{i_1},\ldots,G_{i_{2k}}\}$.
Without loss of generality, we assume that $G_{i_{2k+1}}=\{1,2,3\}$. By (ii), we know that for any $j=i_{2k+2},\ldots,i_{n_3}$,
$$
G_j\in \{\{1,2,4\}, \{1,3,4\}, \{2,3,4\},\{1,2,5\}, \{1,3,5\}, \{2,3,5\},\{1,2,6\}, \{1,3,6\}, \{2,3,6\}\}.
$$
Denote
\begin{eqnarray*}
&&\cal{H}_{4}=\{\{1,2,4\}, \{1,3,4\}, \{2,3,4\}\},\\
&&\cal{H}_5=\{\{1,2,5\}, \{1,3,5\}, \{2,3,5\}\},\\
&&\cal{H}_6=\{\{1,2,6\}, \{1,3,6\}, \{2,3,6\}\}.
\end{eqnarray*}
For any $A\in \cal{H}_4,B\in \cal{H}_5,C\in\cal{H}_6$, we have
\begin{eqnarray*}
&&\{1,2,3\}\cup A\cup B=\{1,2,3,4,5\},\\
 &&\{1,2,3\}\cup A\cup C=\{1,2,3,4,6\},\\
 && \{1,2,3\}\cup B\cup C=\{1,2,3,5,6\}.
\end{eqnarray*}
By $\cal{F}=\{\emptyset, M_6\}\cup \cal{G}_3\cup \cal{G}_4$, without loss of generality, we can assume that
$$
\cal{F}\cap \cal{H}_4\neq \emptyset, \cal{F}\cap \cal{H}_5=\cal{F}\cap \cal{H}_6=\emptyset.
$$
Now by Lemma \ref{lem-1.2}, we know that all the 4 elements in $\{1,2,3,4\}$ belong to at least $n_3-2k-1$ set(s) in $\{G_{i_{2k+1}},\ldots,G_{i_{n_3}}\}$, and thus belong to at least half of the sets in $\cal{G}_3$.
\end{itemize}
\end{itemize}
\end{itemize}

\subsection{$\cal{F}=\{\emptyset, M_6\}\cup \cal{G}_3\cup \cal{G}_4\cup \cal{G}_5$}

Denote $\cal{G}_3=\{G_1,\ldots,G_{n_3}\}, \cal{G}_4=\{H_1,\ldots,H_{n_4}\}$ and $\cal{G}_5=\{I_1,\ldots,I_{n_5}\}$.

\begin{itemize}
\item[(1)] $n_5=1$. Now $\cal{G}_5=\{I_1\}$.  Without loss of generality, we assume that $I_1=\{1,2,3,4,5\}$.

\begin{itemize}
\item[(1.1)] $n_4=1$. Now $\cal{G}_4=\{H_1\}$.

\begin{itemize}
\item[(1.1.1)]  $n_3=1$. Now $\cal{G}_3=\{G_1\}$.  If $H_1\subset I_1$, then all the 4 elements in $H_1$ belong to at least two sets among the three sets in $\cal{G}_3\cup \cal{G}_4\cup \cal{G}_5$ and
thus belong to at least half of the sets in $\cal{F}$. If $H_1\nsubseteq I_1$, then $H_1\cup I_1=M_6$ and thus in this case all the 3 elements in $G_1$ belong to at least half of the sets in $\cal{F}$.

\item[(1.1.2)]  $n_3\geq 2$. Now for any $i,j=1,\ldots,n_3,i\neq j$, we have $G_i\cup G_j=M_6$ or $G_i\cup G_j=I_1=\{1,2,3,4,5\}$ or $G_i\cup G_j=H_1$.
\begin{itemize}
\item[(1.1.2.1)]    $n_3$ is an even number and there exists a permutation $(i_1,\ldots,i_{n_3})$ of $(1,\ldots,n_3)$ such that $G_{i_1}\cup G_{i_2}=\cdots=G_{i_{n_3-1}}\cup G_{i_{n_3}}=M_6$. Then all the 6 elements
in $M_6$ belong to half of the sets in $\cal{G}_3$. Hence all the elements in $H_1\cup I_1$ belong to at least half of the sets in $\cal{F}$.

\item[(1.1.2.2)]    $n_3$ is an odd number and there exists a permutation $(i_1,\ldots,i_{n_3})$ of $(1,\ldots,n_3)$ such that $G_{i_1}\cup G_{i_2}=\cdots=G_{i_{n_3-2}}\cup G_{i_{n_3-1}}=M_6$.
Then all the 6 elements
in $M_6$ belong to half of the sets in $\{G_{i_1},\ldots,G_{i_{n_3-1}}\}$.  If $H_1\subset I_1$, then all the 4 elements in $H_1$ belong to at least two sets among the three sets in
$\{G_{i_{n_3}},H_1,I_1\}$ and thus belong to at least half of the sets in $\cal{F}$.  If $H_1\nsubseteq I_1$, then $H_1\cup I_1=M_6$ and thus in this case all the 3 elements in $G_{i_{n_3}}$
 belong to at least half of the sets in $\cal{F}$.

\item[(1.1.2.3)]  We can decompose $\cal{G}_3$ into two disjoint parts $\{G_{i_1},\ldots,G_{i_{2k}}\}$ and $\{G_{i_{2k+1}},\ldots,G_{i_{n_3}}\}$, where
$\{i_1,\ldots,i_{n_3}\}=\{1,\ldots,n_3\}$, $n_3-2k\geq 2$, and

\quad (i) $G_{i_1}\cup G_{i_2}=\cdots=G_{i_{2k-1}}\cup G_{i_{2k}}=M_6$;

\quad (ii) for any two different indexes $\{i,j\}$ from $\{i_{2k+1},\ldots,i_{n_3}\}, G_i\cup G_j=I_1=\{1,2,3,4,5\}$ or $G_i\cup G_j=H_1$.

Then all the 6 elements in $M_6$ belong to  half of the sets in $\{G_{i_1},\ldots,G_{i_{2k}}\}$.\smallskip

\quad (a)  $H_1\subset \{1,2,3,4,5\}$. Without loss of generality, we assume that $H_1=\{1,2,3,4\}$.

\quad\quad (a.1) $n_3-2k$ is an even number and there exists a permutation $(j_1,\ldots,j_{n_3-2k})$ of $(i_{2k+1},\ldots,i_{n_3})$ such that
$G_{j_1}\cup G_{j_2}=\cdots=G_{j_{n_3-2k-1}}\cup G_{j_{n_3-2k}}=\{1,2,3,4,5\}$. Then all the 5 elements in $\{1,2,3,4,5\}$ belong to at least half of the sets in $\{G_{i_{2k+1}},\ldots,G_{i_{n_3}}\}$.
Note that all the 5 elements in $\{1,2,3,4,5\}$ belong to at least one of the two sets $H_1$ and $I_1$. Then we know that  in this case, all the 5 elements in $\{1,2,3,4,5\}$ belong to at least half of the sets in $\cal{F}$.\smallskip

\quad\quad (a.2)  $n_3-2k$ is an odd number and there exists a permutation $(j_1,\ldots,j_{n_3-2k})$ of $(i_{2k+1},\ldots,i_{n_3})$ such that
$G_{j_1}\cup G_{j_2}=\cdots=G_{j_{n_3-2k-2}}\cup G_{j_{n_3-2k-1}}=\{1,2,3,4,5\}$.  Then all the 5 elements in $\{1,2,3,4,5\}$ belong to at least half of the sets in
$\{G_{j_1},\ldots,G_{j_{n_3-2k-1}}\}$. Note that all the 4 elements in $\{1,2,3,4\}$ belong to at least two sets among the three sets $\{G_{j_{n_3-2k}}, H_1,I_1\}$.
Then we know that in this case, all the 4 elements in $\{1,2,3,4\}$ belong to at least half of the sets in $\cal{F}$.\smallskip

\quad\quad (a.3) We can decompose $\{G_{i_{2k+1}},\ldots,G_{i_{n_3}}\}$ into two disjoint parts $\{G_{j_1},\ldots,\linebreak G_{j_{2l}}\}$ and $\{G_{j_{2l+1}},\ldots,G_{j_{n_3-2k}}\}$, where $\{j_1,\ldots,j_{n_3-2k}\}=
\{i_{2k+1},\ldots,i_{n_3}\}$, $n_3-2k-2l\geq 2$, and

\quad\quad\quad (iii) $G_{j_1}\cup G_{j_2}=\cdots=G_{j_{2l-1}}\cup G_{j_{2l}}=\{1,2,3,4,5\}$;

\quad\quad\quad (iv) for any two different indexes $\{i,j\}$ from $\{j_{2l+1},\ldots,j_{n_3-2k}\}, G_i\cup G_j=H_1=\{1,2,3,4\}$.

Then all the 5 elements in $\{1,2,3,4,5\}$ belong to at least half of the sets in $\{G_{j_1},\ldots,G_{j_{2l}}\}$. By Lemma \ref{lem-1.2}, we know that all the 4 elements in $\{1,2,3,4\}$ belong to at least $n_3-2k-2l-1$ set(s)
in $\{G_{j_{2l+1}},\ldots,G_{j_{n_3-2k}}\}$ and thus belong to at least half of the sets in $\{G_{j_{2l+1}},\ldots,G_{j_{n_3-2k}}\}$. Hence in this case, all the 4 elements in $\{1,2,3,4\}$ belong to at least half of the sets in
$\cal{F}$.
\smallskip

\quad (b)  $H_1\nsubseteq \{1,2,3,4,5\}$. Without loss of generality, we assume that $H_1=\{1,2,3,6\}$. By following the proof in (a), we can get that in this case
all the 3 elements in $\{1,2,3\}$ belong to at least half of the sets in $\cal{F}$.
\end{itemize}
\end{itemize}

\item[(1.2)] $n_4\geq 2$.  Now for any $i,j=1,\ldots,n_4,i\neq j$, we have $H_i\cup H_j=M_6$ or $H_i\cup H_j=I_1=\{1,2,3,4,5\}$.

\begin{itemize}
\item[(1.2.1)]   $n_4$ is an even number and there exists a permutation $(i_1,\ldots,i_{n_4})$ of $(1,\ldots,n_4)$ such that $H_{i_1}\cup H_{i_2}=\cdots=H_{i_{n_4-1}}\cup H_{i_{n_4}}=M_6$. Then all the 6 elements in
$M_6$ belong to at least half of the sets in $\cal{G}_4$.

\begin{itemize}
\item[(1.2.1.1)]  $n_3=1$. Now $\cal{G}_3=\{G_1\}$. In this case, all the elements in $G_1\cup I_1$ belong to at least half of the sets in $\cal{F}$.

\item[(1.2.1.2)]  $n_3\geq 2$. Now for any $i,j=1,\ldots,n_3,i\neq j$, we have $G_i\cup G_j=M_6$ or $G_i\cup G_j=I_1=\{1,2,3,4,5\}$ or $|G_i\cup G_j|=4$.

\quad (a)    $n_3$ is an even number and there exists a permutation $(j_1,\ldots,j_{n_3})$ of $(1,\ldots,n_3)$ such that $G_{j_1}\cup G_{j_2}=\cdots=G_{j_{n_3-1}}\cup G_{j_{n_3}}=M_6$. Then all the 6 elements
in $M_6$ belong to half of the sets in $\cal{G}_3$. Hence in this case all the 5 elements in $I_1$ belong to at least half of the sets in $\cal{F}$.\\

\quad (b)  $n_3$ is an odd number and there exists a permutation $(j_1,\ldots,j_{n_3})$ of $(1,\ldots,n_3)$ such that $G_{j_1}\cup G_{j_2}=\cdots=G_{j_{n_3-2}}\cup G_{j_{n_3-1}}=M_6$.  Then all the 6 elements in $M_6$ belong to half of the sets in $\{G_{j_1},\ldots,G_{j_{n_3-1}}\}$. Hence in this case,  all the  elements in $G_{j_{n_3}}\cup I_1$  belong to at least half of the sets in $\cal{F}$.\\

\quad (c) We can decompose $\cal{G}_3$ into two disjoint parts $\{G_{j_1},\ldots,G_{j_{2k}}\}$ and\linebreak $\{G_{j_{2k+1}},\ldots,G_{j_{n_3}}\}$, where
$\{j_1,\ldots,j_{n_3}\}=\{1,\ldots,n_3\}$, $n_3-2k\geq 2$, and

\quad\quad (i) $G_{j_1}\cup G_{j_2}=\cdots=G_{j_{2k-1}}\cup G_{j_{2k}}=M_6$;

\quad\quad (ii) for any two different indexes $\{i,j\}$ from $\{j_{2k+1},\ldots,j_{n_3}\}, G_i\cup G_j=I_1=\{1,2,3,4,5\}$ or $|G_i\cup G_j|=4$.

Then all the 6 elements in $M_6$ belong to  half of the sets in $\{G_{j_1},\ldots,G_{j_{2k}}\}$.\smallskip

\quad (c.1) $n_3-2k$ is an even number and there exists a permutation $(m_1,\ldots,m_{n_3-2k})$ of $(j_{2k+1},\ldots,j_{n_3})$ such that
$G_{m_1}\cup G_{m_2}=\cdots=G_{m_{n_3-2k-1}}\cup G_{m_{n_3-2k}}=\{1,2,3,4,5\}$. Hence in this case, all the 5 elements in $\{1,2,3,4,5\}$ belong to at least half of the sets in $\cal{F}$.\smallskip

\quad (c.2) $n_3-2k$ is an odd number and there exists a permutation $(m_1,\ldots,m_{n_3-2k})$ of $(j_{2k+1},\ldots,j_{n_3})$ such that
$G_{m_1}\cup G_{m_2}=\cdots=G_{m_{n_3-2k-2}}\cup G_{m_{n_3-2k-1}}=\{1,2,3,4,5\}$.  Then all the 5 elements in $\{1,2,3,4,5\}$ belong to at least half of the sets in
$\{G_{m_1},\ldots,G_{m_{n_3-2k-1}}\}$. Note that all the 5 elements in $\{1,2,3,4,5\}$ belong to at least one of the two sets $G_{m_{n_3-2k}}$ and $\{1,2,3,4,5\}$. Hence in this case, all the 5 elements in $\{1,2,3,4,5\}$ belong to at least half of the sets in $\cal{F}$.
\smallskip

\quad (c.3)  We can decompose $\{G_{j_{2k+1}},\ldots,G_{j_{n_3}}\}$ into two disjoint parts $\{G_{m_1},\ldots,\linebreak G_{m_{2l}}\}$ and $\{G_{m_{2l+1}},\ldots,G_{m_{n_3-2k}}\}$, where $\{m_1,\ldots,m_{n_3-2k}\}=
\{j_{2k+1},\ldots,j_{n_3}\}$, $n_3-2k-2l\geq 2$, and

\quad\quad (iii) $G_{m_1}\cup G_{m_2}=\cdots=G_{m_{2l-1}}\cup G_{m_{2l}}=\{1,2,3,4,5\}$;

\quad\quad (iv) for any two different indexes $\{i,j\}$ from $\{m_{2l+1},\ldots,m_{n_3-2k}\}, G_i\cup G_j\in \cal{G}_4$.

Then all the 5 elements in $\{1,2,3,4,5\}$ belong to at least half of the sets in $\{G_{m_1},\ldots,G_{m_{2l}}\}$.
\smallskip

\quad (c.3.1) There exists $t\in \{m_{2l+1},\ldots,m_{n_3-2k}\}$ such that $G_t\subset \{1,2,3,4,5\}$. Without loss of generality, we assume that
$G_{m_{2l+1}}=\{1,2,3\}$. Then for any $t\in \{m_{2l+2},\ldots,
m_{n_3-2k}\}$, we have $G_t\in \cal{H}_4\cup \cal{H}_5\cup \cal{H}_6$, where
\begin{eqnarray*}
&&\cal{H}_4=\{\{1,2,4\}, \{1,3,4\},\{2,3,4\}\},\\
&&\cal{H}_5=\{\{1,2,5\}, \{1,3,5\},\{2,3,5\}\},\\
&&\cal{H}_6=\{\{1,2,6\}, \{1,3,6\},\{2,3,6\}\}.
\end{eqnarray*}
For simplicity, define $\cal{H}:=\{G_{m_{2l+1}},\ldots,G_{m_{n_3-2k}}\}$.  We have the following 7 cases:
\smallskip

\quad\quad (c.3.1.1)  $\cal{H}\cap \cal{H}_4\neq \emptyset, \cal{H}\cap \cal{H}_5=\cal{H}\cap \cal{H}_6=\emptyset$.

\quad\quad (c.3.1.2)  $\cal{H}\cap \cal{H}_5\neq \emptyset, \cal{H}\cap \cal{H}_4=\cal{H}\cap \cal{H}_6=\emptyset$.

\quad\quad (c.3.1.3)  $\cal{H}\cap \cal{H}_6\neq \emptyset, \cal{H}\cap \cal{H}_4=\cal{H}\cap \cal{H}_5=\emptyset$.

\quad\quad (c.3.1.4)  $\cal{H}\cap \cal{H}_4\neq \emptyset, \cal{H}\cap \cal{H}_5\neq \emptyset, \cal{H}\cap \cal{H}_6=\emptyset$.

\quad\quad (c.3.1.5)  $\cal{H}\cap \cal{H}_4\neq \emptyset, \cal{H}\cap \cal{H}_6\neq \emptyset, \cal{H}\cap \cal{H}_5=\emptyset$.

\quad\quad (c.3.1.6)  $\cal{H}\cap \cal{H}_5\neq \emptyset, \cal{H}\cap \cal{H}_6\neq \emptyset, \cal{H}\cap \cal{H}_4=\emptyset$.

\quad\quad (c.3.1.7)  $\cal{H}\cap \cal{H}_4\neq \emptyset, \cal{H}\cap \cal{H}_5\neq \emptyset, \cal{H}\cap \cal{H}_6\neq \emptyset$.
\smallskip

\quad As to (c.3.1.1), by Lemma \ref{lem-1.2}, we know that all the 4 elements in $\{1,2,3,4\}$ belong to at least $n_3-2k-1$ set(s) in
$\cal{H}$ and thus belong to at least half of the sets in $\cal{H}$. Hence in this case, all the 4 elements in $\{1,2,3,4\}$ belong to at least half of the sets in $\cal{F}$.
\smallskip

\quad As to (c.3.1.2) and  (c.3.1.3), we can get that all the 3 elements in $\{1,2,3\}$  belong to at least half of the sets in $\cal{F}$. \smallskip

\quad As to (c.3.1.4), without loss of generality, we assume  that $\{1,2,4\}\in \cal{H}\cap \cal{H}_4$. Then by (iv), we know that $\cal{H}\cap \cal{H}_5=\{\{1,2,5\}\}$ and by (iv) again we get that
$\cal{H}\cap \cal{H}_4=\{\{1,2,4\}\}$. Thus in this case $\cal{H}=\{\{1,2,3\}, \{1,2,4\}, \{1,2,5\}\}$. Note that all the 5 elements in $\{1,2,3,4,5\}$ belong to at least two sets among the 4 sets in $\{\{1,2,3\}, \{1,2,4\}, \{1,2,5\},
\{1,2,3,4,5\}\}$. Then we obtain that all the 5 elements in $\{1,2,3,4,5\}$ belong to at least half of the sets in $\cal{F}$.
\smallskip

\quad As to (c.3.1.5), without loss of generality, we assume  that $\{1,2,4\}\in \cal{H}\cap \cal{H}_4$. Then by (iv), we know that $\cal{H}\cap \cal{H}_6=\{\{1,2,6\}\}$ and by (iv) again we get that
$\cal{H}\cap \cal{H}_4=\{\{1,2,4\}\}$. Thus in this case $\cal{H}=\{\{1,2,3\}, \{1,2,4\}, \{1,2,6\}\}$. By $\{1,2,3\}\cup \{1,2,4\}\cup \{1,2,6\}=\{1,2,3,4,6\}\in \cal{G}_5=\{\{1,2,3,4,5\}\}$, we know that this case
 is impossible.

\quad As to (c.3.1.6) and  (c.3.1.7), by following the analysis to (c.3.1.5), we know that these  two  cases are impossible. \smallskip

\bigskip

\quad (c.3.2) For any  $t\in  \{m_{2l+1},\ldots,m_{n_3-2k}\}$,  $G_t\nsubseteq \{1,2,3,4,5\}$.  Without loss of generality, we assume that $G_{m_{2l+1}}=\{1,2,6\}$. Then by (iv), we know that for any $t=m_{2l+2},\ldots,m_{n_3-2k}$, we have
$$
G_t\in \{\{1,3,6\}, \{1,4,6\}, \{1,5,6\}, \{2,3,6\} ,\{2,4,6\}, \{2,5,6\}\}.
$$
Without loss of generality, we assume that $\{1,3,6\}\in \{G_{m_{2l+2}},\ldots,G_{m_{n_3-2k}}\}$. Then by $\cal{G}_5=\{\{1,2,3,4,5\}\}$, we need only to consider the following two subcases:

\quad\quad  (c.3.2.1)  $\{G_{m_{2l+1}},\ldots,G_{m_{n_3-2k}}\}=\{\{1,2,6\}, \{1,3,6\}\}$.

\quad\quad (c.3.2.2)  $\{G_{m_{2l+1}},\ldots,G_{m_{n_3-2k}}\}=\{\{1,2,6\}, \{1,3,6\}, \{2,3,6\}\}$.
\smallskip

\quad As to (c.3.2.1), all the 4 elements in $\{1,2,3,6\}$ belong to at least two sets among the 3 sets in $\{\{1,2,6\}, \{1,3,6\}, \{1,2,3,4,5\}\}$. Hence in this case all the 3 elements in $\{1,2,3,6\}\cap \{1,2,3,4,5\}$
(i.e. $\{1,2,3\}$) belong to at least half of the sets in $\cal{F}$. \smallskip

\quad As to (c.3.2.2), we can easily know that all the 3 elements in $\{1,2,3\}$ belong to at least half of the sets in $\cal{F}$.
\end{itemize}

\item[(1.2.2)]   $n_4$ is an odd number and there exists a permutation $(i_1,\ldots,i_{n_4})$ of $(1,\ldots,n_4)$ such that $H_{i_1}\cup H_{i_2}=\cdots=H_{i_{n_4-2}}\cup H_{i_{n_4-1}}=M_6$.
Then all the 6 elements in $M_6$ belong to at least half of the sets in $\{H_{i_1},\ldots,H_{i_{n_4-1}}\}$.

\begin{itemize}
\item[(1.2.2.1)]  $n_3=1$. Now $\cal{G}_3=\{G_1\}$. If $H_{i_{n_4}}\subset I_1$, then all the 4 elements in $H_{i_{n_4}}$ belong to at least two sets among the
three sets in $\{G_1,H_{i_{n_4}},I_1\}$ and thus all the 4 elements in $H_{i_{n_4}}$ belong to at least half of the sets in $\cal{F}$. If $H_{i_{n_4}}\nsubseteq I_1$, then $H_{i_{n_4}}\cup I_1=M_6$ and thus in this
case all the 3 elements in $G_1$ belong to at least half of the sets in $\cal{F}$.

\item[(1.2.2.2)]  $n_3\geq 2$.  Now for any $i,j=1,\ldots,n_3,i\neq j$, we have $G_i\cup G_j=M_6$ or $G_i\cup G_j=I_1=\{1,2,3,4,5\}$ or $|G_i\cup G_j|=4$.

\quad (a)    $n_3$ is an even number and there exists a permutation $(j_1,\ldots,j_{n_3})$ of $(1,\ldots,n_3)$ such that $G_{j_1}\cup G_{j_2}=\cdots=G_{j_{n_3-1}}\cup G_{j_{n_3}}=M_6$. Then all the 6 elements
in $M_6$ belong to half of the sets in $\cal{G}_3$. Hence in this case all the  elements in $H_{i_{n_4}}\cup I_1$ belong to at least half of the sets in $\cal{F}$.

\bigskip

\quad (b)  $n_3$ is an odd number and there exists a permutation $(j_1,\ldots,j_{n_3})$ of $(1,\ldots,n_3)$ such that $G_{j_1}\cup G_{j_2}=\cdots=G_{j_{n_3-2}}\cup G_{j_{n_3-1}}=M_6$. Then all the 6 elements
in $M_6$ belong to half of the sets in $\{G_{j_1},\ldots,G_{j_{n_3-1}}\}$. Now if $H_{i_{n_4}}\subset I_1$, then all the 4 elements in $H_{i_{n_4}}$ belong to at least two sets among the
three sets in $\{G_{j_{n_3}},H_{i_{n_4}},I_1\}$ and thus all the 4 elements in $H_{i_{n_4}}$ belong to at least half of the sets in $\cal{F}$. If $H_{i_{n_4}}\nsubseteq I_1$, then $H_{i_{n_4}}\cup I_1=M_6$ and thus in this  case all the 3 elements in $G_{j_{n_3}}$ belong to at least half of the sets in $\cal{F}$.
\bigskip

\quad (c) We can decompose $\cal{G}_3$ into two disjoint parts $\{G_{j_1},\ldots,G_{j_{2k}}\}$ and\linebreak $\{G_{j_{2k+1}},\ldots,G_{j_{n_3}}\}$, where
$\{j_1,\ldots,j_{n_3}\}=\{1,\ldots,n_3\}$, $n_3-2k\geq 2$, and

\quad\quad (i) $G_{j_1}\cup G_{j_2}=\cdots=G_{j_{2k-1}}\cup G_{j_{2k}}=M_6$;

\quad\quad (ii) for any two different indexes $\{i,j\}$ from $\{j_{2k+1},\ldots,j_{n_3}\}, G_i\cup G_j=I_1=\{1,2,3,4,5\}$ or $|G_i\cup G_j|=4$.

Then all the 6 elements in $M_6$ belong to  half of the sets in $\{G_{j_1},\ldots,G_{j_{2k}}\}$.
\bigskip

\quad (c.1) $n_3-2k$ is an even number and there exists a permutation $(m_1,\ldots,m_{n_3-2k})$ of $(j_{2k+1},\ldots,j_{n_3})$ such that
$G_{m_1}\cup G_{m_2}=\cdots=G_{m_{n_3-2k-1}}\cup G_{m_{n_3-2k}}=\{1,2,3,4,5\}$, which together with the fact that all the
5 elements in $\{1,2,3,4,5\}$ belong to at least one of the two sets $H_{i_{n_4}}$ and $\{1,2,3,4,5\}$, implies that  all the 5 elements in $\{1,2,3,4,5\}$ belong to at least half of the sets in $\cal{F}$.
\bigskip

\quad (c.2) $n_3-2k$ is an odd number and there exists a permutation $(m_1,\ldots,m_{n_3-2k})$ of $(j_{2k+1},\ldots,j_{n_3})$ such that
$G_{m_1}\cup G_{m_2}=\cdots=G_{m_{n_3-2k-2}}\cup G_{m_{n_3-2k-1}}=\{1,2,3,4,5\}$.  Then all the 5 elements in $\{1,2,3,4,5\}$ belong to at least half of the sets in
$\{G_{m_1},\ldots,G_{m_{n_3-2k-1}}\}$.
\smallskip

\quad (c.2.1)  If $H_{i_{n_4}}\subset \{1,2,3,4,5\}$, then all the 4 elements in $H_{i_{n_4}}$ belong to at least two sets among the three sets in
$\{G_{m_{n_3-2k}},H_{i_{n_4}},\{1,2,3,4,5\}\}$. Hence in this case, all the 4 elements in $H_{i_{n_4}}$ belong to at least half of the sets in $\cal{F}$.
\bigskip

\quad (c.2.2)   If $H_{i_{n_4}}\nsubseteq \{1,2,3,4,5\}$, then  $H_{i_{n_4}}\cup \{1,2,3,4,5\}=M_6$. Without loss of generality, we assume that $H_{i_{n_4}}=\{1,2,3,6\}$.
\bigskip

\quad  (c.2.2.1)  $G_{m_{n_3-2k}}\subset \{1,2,3,4,5\}$.   Then all the 3 elements in $G_{m_{n_3-2k}}$  belong to at least two sets among the three sets in $\{G_{m_{n_3-2k}},H_{i_{n_4}},\{1,2,3,4,5\}\}$ and thus belong to at least half of the sets in $\cal{F}$.
\bigskip

\quad (c.2.2.2)   $G_{m_{n_3-2k}}\nsubseteq \{1,2,3,4,5\}$. Then $G_{m_{n_3-2k}}\cup \{1,2,3,4,5\}=M_6$.  Hence all the 6 elements in $M_6$ belong to at least one of the two sets $G_{m_{n_3-2k}}$ and $\{1,2,3,4,5\}$ and thus all the 4 elements in $H_{i_{n_4}}$ belong to at least two sets among the three sets in $\{G_{m_{n_3-2k}},H_{i_{n_4}},\{1,2,3,4,5\}\}$.  Hence in this case all the elements in
$H_{i_{n_4}}\cap \{1,2,3,4,5\}$ (i.e. $\{1,2,3\}$) belong to at least half of the sets in $\cal{F}$.
\bigskip

\quad (c.3)  We can decompose $\{G_{j_{2k+1}},\ldots,G_{j_{n_3}}\}$ into two disjoint parts\linebreak
 $\{G_{m_1},\ldots,G_{m_{2l}}\}$ and $\{G_{m_{2l+1}},\ldots,G_{m_{n_3-2k}}\}$, where $\{m_1,\ldots,m_{n_3-2k}\}=
\{j_{2k+1},\linebreak \ldots, j_{n_3}\}$, $n_3-2k-2l\geq 2$, and

\quad\quad (iii) $G_{m_1}\cup G_{m_2}=\cdots=G_{m_{2l-1}}\cup G_{m_{2l}}=\{1,2,3,4,5\}$;

\quad\quad (iv) for any two different indexes $\{i,j\}$ from $\{m_{2l+1},\ldots,m_{n_3-2k}\}, G_i\cup G_j\in \cal{G}_4$.

Then all the 5 elements in $\{1,2,3,4,5\}$ belong to at least half of the sets in $\{G_{m_1},\ldots,G_{m_{2l}}\}$.
\smallskip

\quad (c.3.1) $H_{i_{n_4}}\subset \{1,2,3,4,5\}$. 

\quad  (c.3.1.1) There exists $t\in \{m_{2l+1},\ldots,m_{n_3-2k}\}$ such that $G_t\subset \{1,2,3,4,5\}$. Without loss of generality, we assume that $G_{m_{2l+1}}=\{1,2,3\}$.
Then for any $t\in \{m_{2l+2},\ldots, m_{n_3-2k}\}$, we have $G_t\in \cal{H}_4\cup \cal{H}_5\cup \cal{H}_6$, where
\begin{eqnarray*}
&&\cal{H}_4=\{\{1,2,4\}, \{1,3,4\},\{2,3,4\}\},\\
&&\cal{H}_5=\{\{1,2,5\}, \{1,3,5\},\{2,3,5\}\},\\
&&\cal{H}_6=\{\{1,2,6\}, \{1,3,6\},\{2,3,6\}\}.
\end{eqnarray*}
For simplicity, define $\cal{H}:=\{G_{m_{2l+1}},\ldots,G_{m_{n_3-2k}}\}$.  We have the following 7 cases:
\smallskip

\quad\quad (c.3.1.1.1)  $\cal{H}\cap \cal{H}_4\neq \emptyset, \cal{H}\cap \cal{H}_5=\cal{H}\cap \cal{H}_6=\emptyset$.

\quad\quad (c.3.1.1.2)  $\cal{H}\cap \cal{H}_5\neq \emptyset, \cal{H}\cap \cal{H}_4=\cal{H}\cap \cal{H}_6=\emptyset$.

\quad\quad (c.3.1.1.3)  $\cal{H}\cap \cal{H}_6\neq \emptyset, \cal{H}\cap \cal{H}_4=\cal{H}\cap \cal{H}_5=\emptyset$.

\quad\quad (c.3.1.1.4)  $\cal{H}\cap \cal{H}_4\neq \emptyset, \cal{H}\cap \cal{H}_5\neq \emptyset, \cal{H}\cap \cal{H}_6=\emptyset$.

\quad\quad (c.3.1.1.5)  $\cal{H}\cap \cal{H}_4\neq \emptyset, \cal{H}\cap \cal{H}_6\neq \emptyset, \cal{H}\cap \cal{H}_5=\emptyset$.

\quad\quad (c.3.1.1.6)  $\cal{H}\cap \cal{H}_5\neq \emptyset, \cal{H}\cap \cal{H}_6\neq \emptyset, \cal{H}\cap \cal{H}_4=\emptyset$.

\quad\quad (c.3.1.1.7)  $\cal{H}\cap \cal{H}_4\neq \emptyset, \cal{H}\cap \cal{H}_5\neq \emptyset, \cal{H}\cap \cal{H}_6\neq \emptyset$.

By the analysis in (1.2.1.2)(c.3.1), we need only to consider the first four cases (c.3.1.1.1)-(c.3.1.1.4).\smallskip

\quad As to (c.3.1.1.1), by Lemma \ref{lem-1.2}, we know that all the 4 elements in $\{1,2,3,4\}$ belong to at least half of the sets in $\cal{H}$. Since $H_{i_{n_4}}\subset \{1,2,3,4,5\}$,
we know that $|H_{i_{n_4}}\cap \{1,2,3,4\}|\geq 3$. Hence in this case, all the elements in $H_{i_{n_4}}\cap \{1,2,3,4\}$ belong to at least half of the sets in $\cal{F}$.
\smallskip

\quad As to (c.3.1.1.2) ,  by following the analysis in (c.3.1.1.1), we get that $|H_{i_{n_4}}\cap \{1,2,3,5\}|\geq 3$ and all the elements in $H_{i_{n_4}}\cap \{1,2,3,5\}$ belong to at least half of the sets in $\cal{F}$.
\smallskip

\quad As to (c.3.1.1.3),  by Lemma \ref{lem-1.2}, we know that all the 4 elements in $\{1,2,3,6\}$ belong to at least half of the sets in $\cal{H}$. We have the following 3 subcases:

\quad\quad  (c.3.1.1.3-1)  $\cal{H}\cap \cal{H}_6=1$. Take $\cal{H}\cap \cal{H}_6=\{\{1,2,6\}\}$ for example.  By $H_{i_{n_4}}\subset \{1,2,3,4,5\}$, we know that
$$
H_{i_{n_4}}\in \{\{1,2,3,4\},\{1,2,3,5\}, \{1,2,4,5\},\{1,3,4,5\}, \{2,3,4,5\}\}.
$$

\quad If $H_{i_{n_4}}=\{1,2,3,4\}$, then $\{1,2,6\}\cup \{1,2,3,4\}=\{1,2,3,4,6\}\in\cal{G}_5=\{\{1,2,3,4,5\}\}$. It is impossible. Similarly, if $H_{i_{n_4}}\in \{\{1,2,3,5\},\{1,2,4,5\}\}$, it is impossible.

\quad If $H_{i_{n_4}}\in \{\{1,3,4,5\}, \{2,3,4,5\}\},$ then $H_{i_{n_4}}\cup \{1,2,6\}=M_6$.  Now all the 3 elements in $\{1,2,3\}$ belong to at least half of the sets in $\cal{F}$.
\smallskip

\quad\quad  (c.3.1.1.3-2)  $\cal{H}\cap \cal{H}_6=2$. Take $\cal{H}\cap \cal{H}_6=\{\{1,2,6\},\{1,3,6\}\}$ for example. By following analysis in (c.3.1.1.3-1),
 we get that all the 3 elements in $\{1,2,3\}$ belong to at least half of the sets in $\cal{F}$.
\smallskip

\quad\quad  (c.3.1.1.3-2)  $\cal{H}\cap \cal{H}_6=3$. Now $\cal{H}\cap \cal{H}_6=\{\{1,2,6\},\{1,3,6\}, \{2,3,6\}\}$. In this case, by $\cal{G}_5=\{\{1,2,3,4,5\}\}$, we must have $H_{i_{n_4}}=\{1,3,4,5\}$. It is easy to check that in this case
all the 3 elements in $\{1,2,3\}$ belong to at least half of the sets in $\cal{F}$.

\bigskip

\quad  (c.3.1.2) For any  $t\in  \{m_{2l+1},\ldots,m_{n_3-2k}\}$,  $G_t\nsubseteq \{1,2,3,4,5\}$.  Without loss of generality, we assume that $G_{m_{2l+1}}=\{1,2,6\}$. Then by (iv), we know that for any $t=m_{2l+2},\ldots,m_{n_3-2k}$, we have
$$
G_t\in \{\{1,3,6\}, \{1,4,6\}, \{1,5,6\}, \{2,3,6\} ,\{2,4,6\}, \{2,5,6\}\}.
$$
Without loss of generality, we assume that $\{1,3,6\}\in \{G_{m_{2l+2}},\ldots,G_{m_{n_3-2k}}\}$. Then by $\cal{G}_5=\{\{1,2,3,4,5\}\}$, we need only to consider the following two subcases:

\quad\quad  (c.3.1.1.2-1) $\{G_{m_{2l+1}},\ldots,G_{m_{n_3-2k}}\}=\{\{1,2,6\}, \{1,3,6\}\}$.

\quad\quad (c.3.1.1.2-2) $\{G_{m_{2l+1}},\ldots,G_{m_{n_3-2k}}\}=\{\{1,2,6\}, \{1,3,6\}, \{2,3,6\}\}$.
\smallskip

\quad As to (c.3.1.1.2-1) ,  we know that all the 3 elements in $\{1,2,3\}$ belong to at least 2 sets among the 4 sets in $\{\{1,2,6\},\{1,3,6\}, H_{i_{n_4}}, I_1\}$. Hence in this case,
all the 3 elements in $\{1,2,3\}$ belong to at least half of the sets in $\cal{F}$.  \smallskip

\quad As to (c.3.1.1.2-2) , we can easily know that all the 3 elements in $\{1,2,3\}$ belong to at least 3 sets among the 5 sets  in $\{\{1,2,6\}, \{1,3,6\}, \linebreak
\{2,3,6\},H_{i_{n_4}}, I_1\}$.
Hence in this case,
all the 3 elements in $\{1,2,3\}$ belong to at least half of the sets in $\cal{F}$.  \smallskip

\bigskip

\quad\quad\quad (c.3.2) $H_{i_{n_4}}\nsubseteq \{1,2,3,4,5\}$. Then $H_{i_{n_4}}\cup \{1,2,3,4,5\}=M_6$. By following the analysis in (c.3.1), we can show that there exist 3 elements in $M_6$ which belong to at least half of the sets in $\cal{F}$. We omit the details.
\end{itemize}
\smallskip

\item[(1.2.3)]  We can decompose $\cal{G}_4$ into two disjoint parts $\{H_{i_1},\ldots,H_{i_{2k}}\}$ and $\{H_{i_{2k+1}},\ldots,H_{i_{n_4}}\}$, where $\{i_1,\ldots,i_{n_4}\}=\{1,\ldots,n_4\}$,
$n_4-2k\geq 2$, and

\quad (i) $H_{i_1}\cup H_{i_2}=\cdots=H_{i_{2k-1}}\cup H_{i_{2k}}=M_6$;

\quad (ii) for any two different indexes $\{i,j\}$ from $\{i_{2k+1},\ldots,i_{n_4}\}$, we have $H_i\cup H_j=I_1=\{1,2,3,4,5\}$.

Then all the 6 elements in $M_6$ belong to at least half of the sets in $\{H_{i_1},\ldots,H_{i_{2k}}\}$. By Lemma \ref{lem-1.2}, we know that all the 5 elements in $\{1,2,3,4,5\}$ belong to at least
$n_3-2k-1$ set(s) in $\{H_{i_{2k+1}},\ldots,H_{i_{n_4}}\}$ and thus belong to at least half of the sets in $\{H_{i_{2k+1}},\ldots,H_{i_{n_4}}\}$.

\begin{itemize}
\item[(1.2.3.1)]  $n_3=1$.  Now $\cal{G}_3=\{G_1\}$. Note that all the 5 elements in $I_1=\{1,2,3,4,5\}$ belong to at least one of the two sets $G_1$ and $I_1$. Then we obtain that  all the 5 elements
in $\{1,2,3,4,5\}$ belong to at least half of the sets in $\cal{F}$.

\item[(1.2.3.2)]  $n_3\geq 2$.  Now for any $i,j=1,\ldots,n_3,i\neq j$, we have $G_i\cup G_j=M_6$ or $G_i\cup G_j=I_1=\{1,2,3,4,5\}$ or $|G_i\cup G_j|=4$.

\quad (a)    $n_3$ is an even number and there exists a permutation $(j_1,\ldots,j_{n_3})$ of $(1,\ldots,n_3)$ such that $G_{j_1}\cup G_{j_2}=\cdots=G_{j_{n_3-1}}\cup G_{j_{n_3}}=M_6$. Then all the 6 elements in $M_6$ belong to half of the sets in $\cal{G}_3$. Hence in this case all the 5 elements
in $\{1,2,3,4,5\}$ belong to at least half of the sets in $\cal{F}$.

\bigskip

\quad (b)  $n_3$ is an odd number and there exists a permutation $(j_1,\ldots,j_{n_3})$ of $(1,\ldots,n_3)$ such that $G_{j_1}\cup G_{j_2}=\cdots=G_{j_{n_3-2}}\cup G_{j_{n_3-1}}=M_6$.
Then all the 6 elements in $M_6$ belong to half of the sets in $\{G_{j_1},\ldots,G_{j_{n_3-1}}\}$. Note that all the 5 elements in $I_1=\{1,2,3,4,5\}$ belong to at least one of the two sets
$G_{j_{n_3}}$ and $I_1$. Then we obtain that  all the 5 elements  in $\{1,2,3,4,5\}$ belong to at least half of the sets in $\cal{F}$.
\bigskip

\quad (c) We can decompose $\cal{G}_3$ into two disjoint parts $\{G_{j_1},\ldots,G_{j_{2k}}\}$ and\linebreak $\{G_{j_{2k+1}},\ldots,G_{j_{n_3}}\}$, where
$\{j_1,\ldots,j_{n_3}\}=\{1,\ldots,n_3\}$, $n_3-2k\geq 2$, and

\quad\quad (i) $G_{j_1}\cup G_{j_2}=\cdots=G_{j_{2k-1}}\cup G_{j_{2k}}=M_6$;

\quad\quad (ii) for any two different indexes $\{i,j\}$ from $\{j_{2k+1},\ldots,j_{n_3}\}, G_i\cup G_j=I_1=\{1,2,3,4,5\}$ or $|G_i\cup G_j|=4$.

Then all the 6 elements in $M_6$ belong to  half of the sets in $\{G_{j_1},\ldots,G_{j_{2k}}\}$.
\smallskip

\quad (c.1) $n_3-2k$ is an even number and there exists a permutation $(m_1,\ldots,m_{n_3-2k})$ of $(j_{2k+1},\ldots,j_{n_3})$ such that
$G_{m_1}\cup G_{m_2}=\cdots=G_{m_{n_3-2k-1}}\cup G_{m_{n_3-2k}}=\{1,2,3,4,5\}$. In this case, we get that  all the 5 elements  in $\{1,2,3,4,5\}$ belong to at least half of the sets in $\cal{F}$.
\bigskip

\quad (c.2) $n_3-2k$ is an odd number and there exists a permutation $(m_1,\ldots,m_{n_3-2k})$ of $(j_{2k+1},\ldots,j_{n_3})$ such that
$G_{m_1}\cup G_{m_2}=\cdots=G_{m_{n_3-2k-2}}\cup G_{m_{n_3-2k-1}}=\{1,2,3,4,5\}$.  Then all the 5 elements in $\{1,2,3,4,5\}$ belong to at least half of the sets in
$\{G_{m_1},\ldots,G_{m_{n_3-2k-1}}\}$. Note that all the 5 elements in $I_1=\{1,2,3,4,5\}$ belong to at least one of the two sets
$G_{m_{n_3-2k}}$ and $I_1$. Then we obtain that  all the 5 elements  in $\{1,2,3,4,5\}$ belong to at least half of the sets in $\cal{F}$.
\bigskip

\quad (c.3)  We can decompose $\{G_{j_{2k+1}},\ldots,G_{j_{n_3}}\}$ into two disjoint parts\linebreak
 $\{G_{m_1},\ldots,G_{m_{2l}}\}$ and $\{G_{m_{2l+1}},\ldots,G_{m_{n_3-2k}}\}$, where $\{m_1,\ldots,m_{n_3-2k}\}=
\{j_{2k+1},\linebreak \ldots,j_{n_3}\}$, $n_3-2k-2l\geq 2$, and

\quad\quad (iii) $G_{m_1}\cup G_{m_2}=\cdots=G_{m_{2l-1}}\cup G_{m_{2l}}=\{1,2,3,4,5\}$;

\quad\quad (iv) for any two different indexes $\{i,j\}$ from $\{m_{2l+1},\ldots,m_{n_3-2k}\}, G_i\cup G_j\in \cal{G}_4$.

Then all the 5 elements in $\{1,2,3,4,5\}$ belong to at least half of the sets in $\{G_{m_1},\ldots,G_{m_{2l}}\}$.  By following the analysis in (1.2.1.2)(c.3), we can get that there exist 3
elements in $M_6$ which belong to at least half of the sets in $\cal{F}$.
\end{itemize}
\end{itemize}
\end{itemize}

\item[(2)] $n_5\geq 2$. By Lemma \ref{lem-1.2}, we know that all the 6 elements in $M_6$ belong to at least $n_5-1$ set(s) in $\cal{G}_5$ and thus belong to at least half of
the sets in $\cal{G}_5$. Hence it is enough to show that there exists at least 3 elements in $M_6$ which belong to at least half of the sets in $\cal{G}_3\cup \cal{G}_4$.

\begin{itemize}
\item[(2.1)]  $n_4=1$.  Now $\cal{G}_4=\{H_1\}$.  Without loss of generality, we assume that $H_1=\{1,2,3,4\}$.

\begin{itemize}
\item[(2.1.1)]  $n_3=1$.  Now $\cal{G}_3=\{G_1\}$.  In this case, all the elements in $G_1\cup H_1$ belong to at least half of the sets in $\cal{G}_3\cup \cal{G}_4$.

\item[(2.1.2)]  $n_3\geq 2$. Now for any $i,j=1,\ldots,n_3,i\neq j$, we have $G_i\cup G_j=M_6$ or $G_i\cup G_j\in \cal{G}_5$ or $G_i\cup G_j=H_1$.
\begin{itemize}
\item[(2.1.2.1)]    $n_3$ is an even number and there exists a permutation $(i_1,\ldots,i_{n_3})$ of $(1,\ldots,n_3)$ such that $G_{i_1}\cup G_{i_2}=\cdots=G_{i_{n_3-1}}\cup G_{i_{n_3}}=M_6$. Then all the 6 elements
in $M_6$ belong to half of the sets in $\cal{G}_3$. Hence in this case all the 4 elements in $I_1$ belong to at least half of the sets in $\cal{G}_3\cup \cal{G}_4$.

\item[(2.1.2.2)]    $n_3$ is an odd number and there exists a permutation $(i_1,\ldots,i_{n_3})$ of $(1,\ldots,n_3)$ such that $G_{i_1}\cup G_{i_2}=\cdots=G_{i_{n_3-2}}\cup G_{i_{n_3-1}}=M_6$.
Then all the 6 elements in $M_6$ belong to half of the sets in $\{G_{i_1},\ldots,G_{i_{n_3-1}}\}$.  In this case, all the elements in $G_{i_3}\cup I_1$ belong to at least half of the sets in $\cal{G}_3\cup \cal{G}_4$.

\item[(2.1.2.3)]  We can decompose $\cal{G}_3$ into two disjoint parts $\{G_{i_1},\ldots,\linebreak G_{i_{2k}}\}$ and $\{G_{i_{2k+1}},\ldots,G_{i_{n_3}}\}$, where
$\{i_1,\ldots,i_{n_3}\}=\{1,\ldots,n_3\}$, $n_3-2k\geq 2$, and

\quad (i) $G_{i_1}\cup G_{i_2}=\cdots=G_{i_{2k-1}}\cup G_{i_{2k}}=M_6$;

\quad (ii) for any two different indexes $\{i,j\}$ from $\{i_{2k+1},\ldots,i_{n_3}\}, G_i\cup G_j=H_1=\{1,2,3,4\}$ or $G_i\cup G_j\in \cal{G}_5$.

Then all the 6 elements in $M_6$ belong to  half of the sets in $\{G_{i_1},\ldots,G_{i_{2k}}\}$.
\smallskip

\quad (a) $n_3-2k$ is an even number and there exists a permutation $(j_1,\ldots,j_{n_3-2k})$ of $(i_{2k+1},\ldots,i_{n_3})$ such that
$G_{j_1}\cup G_{j_2}=\cdots=G_{j_{n_3-2k-1}}\cup G_{j_{n_3-2k}}=H_1=\{1,2,3,4\}$. Then all the 4 elements in $\{1,2,3,4\}$ belong to at least half of the sets in $\{G_{i_{2k+1}},\ldots,G_{i_{n_3}}\}$
and thus all the 4 elements in $\{1,2,3,4\}$ belong to at least half of the sets in  $\cal{G}_3\cup \cal{G}_4$.
\smallskip

\quad (b) $n_3-2k$ is an odd number and there exists a permutation $(j_1,\ldots,j_{n_3-2k})$ of $(i_{2k+1},\ldots,i_{n_3})$ such that
$G_{j_1}\cup G_{j_2}=\cdots=G_{j_{n_3-2k-2}}\cup G_{j_{n_3-2k-1}}=H_1=\{1,2,3,4\}$. Then all the 4 elements in $\{1,2,3,4\}$ belong to at least half of the sets in $\{G_{j_1},\ldots,G_{j_{n_3-2k-1}}\}$.
Note that all the  the 4 elements in $\{1,2,3,4\}$ belong to at least one of the two sets $G_{j_{n_3-2k}}$ and $\{1,2,3,4\}$. Then we get that all the 4 elements in $\{1,2,3,4\}$ belong to at least half
 of the sets in  $\cal{G}_3\cup \cal{G}_4$.
 \smallskip

 \quad (c) We can decompose $\{G_{i_{2k+1}},\ldots,G_{i_{n_3}}\}$ into two disjoint two parts $\{G_{j_1},\ldots,\linebreak G_{j_{2l}}\}$ and $\{G_{j_{2l+1}},\ldots,G_{j_{n_3-2k}}\}$, where $\{j_1,\ldots,j_{n_3-2k}\}=\{i_{2k+1},
 \ldots,i_{n_3}\}$, $n_3-2k-2l\geq 2$, and

 \quad\quad  (iii) $G_{j_1}\cup G_{j_2}=\cdots=G_{j_{2l-1}}\cup G_{j_{2l}}=\{1,2,3,4\}$;

 \quad\quad  (iv)  for any two different indexes $\{i,j\}$ from $\{j_{2l+1},\ldots,j_{n_3-2k}\}$, we have $G_i\cup G_j\in \cal{G}_5$.

 Then all the 4 elements in $\{1,2,3,4\}$ belong to at least half of the sets in $\{G_{j_1},\ldots,G_{j_{2l}}\}$. For simplicity, define $\cal{H}:=\{G_{j_{2l+1}},\ldots,G_{j_{n_3-2k}}\}$.   We have the following two subcases:
 \smallskip

 \quad (c.1) There exists $m\in  \{j_{2l+1},\ldots,j_{n_3-2k}\}$ such that $G_m\subset \{1,2,3,4\}$. Without loss of generality, we assume that $G_{j_{2l+1}}=\{1,2,3\}$.
 Then for any $m=j_{2l+2},\ldots,j_{n_3-2k}$, we have  $G_m\in \cal{H}_{4,5}\cup \cal{H}_{4,6}\cup \cal{H}_{5,6}$,
 where
 \begin{eqnarray*}
 &&\cal{H}_{4,5}=\{\{1,4,5\}, \{2,4,5\}, \{3,4,5\}\},\\
 &&\cal{H}_{4,6}=\{\{1,4,6\}, \{2,4,6\}, \{3,4,6\}\},\\
 &&\cal{H}_{5,6}=\{\{1,5,6\}, \{2,5,6\}, \{3,5,6\}\}.
 \end{eqnarray*}
 Since $|\{1,4,5\}\cup \{2,4,5\}|=|\{1,4,5\}\cup \{3,4,5\}|=|\{2,4,5\}\cup \{3,4,5\}|=4$, we get that
 $ |\cal{H}\cap \cal{H}_{4,5}|\leq 1.$  Similarly, we have
$ |\cal{H}\cap \cal{H}_{4,6}|\leq 1,\ \  |\cal{H}\cap \cal{H}_{5,6}|\leq 1.$
Hence we need only to consider the following 7 cases:\smallskip

 \quad\quad (c.1.1) $|\cal{H}\cap \cal{H}_{4,5}|=1, |\cal{H}\cap \cal{H}_{4,6}|=|\cal{H}\cap \cal{H}_{5,6}|=0$.

  \quad\quad (c.1.2) $|\cal{H}\cap \cal{H}_{4,6}|=1, |\cal{H}\cap \cal{H}_{4,5}|=|\cal{H}\cap \cal{H}_{5,6}|=0$.

   \quad\quad (c.1.3) $|\cal{H}\cap \cal{H}_{5,6}|=1, |\cal{H}\cap \cal{H}_{4,5}|=|\cal{H}\cap \cal{H}_{4,6}|=0$.

    \quad\quad (c.1.4) $|\cal{H}\cap \cal{H}_{4,5}|= |\cal{H}\cap \cal{H}_{4,6}|=1, |\cal{H}\cap \cal{H}_{5,6}|=0$.

  \quad\quad (c.1.5) $|\cal{H}\cap \cal{H}_{4,5}|= |\cal{H}\cap \cal{H}_{5,6}|=1, |\cal{H}\cap \cal{H}_{4,6}|=0$.

   \quad\quad (c.1.6) $|\cal{H}\cap \cal{H}_{4,6}|= |\cal{H}\cap \cal{H}_{5,6}|=1, |\cal{H}\cap \cal{H}_{4,5}|=0$.

  \quad\quad (c.1.7) $|\cal{H}\cap \cal{H}_{4,5}|= |\cal{H}\cap \cal{H}_{4,6}|=|\cal{H}\cap \cal{H}_{5,6}|=1$.
  \smallskip

  \quad As to (c.1.1), take $\cal{H}\cap \cal{H}_{4,5}=\{\{1,4,5\}\}$ for example. Now $\cal{H}=\{\{1,2,3\},\linebreak \{1,4,5\}\}$.
   Note that all the 4 elements in $\{1,2,3,4\}$ belong to at least
  two sets among the three sets in $\{\{1,2,3\},\{1,4,5\}, \{1,2,3,4\}\}$. Then we get that in this case, all the 4 elements in $\{1,2,3,4\}$ belong to at least half of the sets
  in $\cal{G}_3\cup \cal{G}_4$.
  \smallskip

  \quad As to (c.1.2), by following the analysis to (c.1.1), we get that all the 4 elements in $\{1,2,3,4\}$ belong to at least half of the sets
  in $\cal{G}_3\cup \cal{G}_4$.
  \smallskip

   \quad As to (c.1.3), by following the analysis to (c.1.1), we get that all the 3 elements in $\{1,2,3\}$ belong to at least half of the sets
  in $\cal{G}_3\cup \cal{G}_4$.
  \smallskip

  \quad As to (c.1.4), take $\cal{H}\cap (\cal{H}_{4,5}\cup \cal{H}_{4,6})=\{\{1,4,5\},\{2,4,6\}\}$ for example. Now $\cal{H}=\{\{1,2,3\}, \{1,4,5\},\{2,4,6\}\}$.
   Note that all the 4 elements in $\{1,2,3,4\}$ belong to at least
  two sets among the four sets in $\{\{1,2,3\},\{1,4,5\}, \{2,4,6\},\linebreak  \{1,2,3,4\}\}$. Then we get that in this case, all the 4 elements in $\{1,2,3,4\}$ belong to at least half of the sets
  in $\cal{G}_3\cup \cal{G}_4$.
  \smallskip

\quad As to (c.1.5), take $\cal{H}\cap (\cal{H}_{4,5}\cup \cal{H}_{4,6})=\{\{1,4,5\},\{2,5,6\}\}$ for example. Now $\cal{H}=\{\{1,2,3\}, \{1,4,5\},\{2,5,6\}\}$.
   Note that all the 4 elements in $\{1,2,3,4\}$ belong to at least
  two sets among the four sets in $\{\{1,2,3\},\{1,4,5\}, \{2,5,6\},\linebreak
   \{1,2,3,4\}\}$. Then we get that in this case, all the 4 elements in $\{1,2,3,4\}$ belong to at least half of the sets
  in $\cal{G}_3\cup \cal{G}_4$.
  \smallskip

  \quad As to (c.1.6), take $\cal{H}\cap (\cal{H}_{4,6}\cup \cal{H}_{5,6})=\{\{1,4,6\},\{2,5,6\}\}$ for example. Now $\cal{H}=\{\{1,2,3\}, \{1,4,6\},\{2,5,6\}\}$.
   Note that all the 4 elements in $\{1,2,3,4\}$ belong to at least
  two sets among the four sets in $\{\{1,2,3\},\{1,4,6\}, \{2,5,6\},\linebreak
   \{1,2,3,4\}\}$. Then we get that in this case, all the 4 elements in $\{1,2,3,4\}$ belong to at least half of the sets
  in $\cal{G}_3\cup \cal{G}_4$.
  \smallskip

 \quad As to (c.1.7), take $\cal{H}\cap \cal{H}_{4,5}=\{\{1,4,5\},\{2,4,6\},\{3,5,6\} \}$ for example. Now $\cal{H}=\{\{1,2,3\}, \{1,4,5\},\{2,4,6\},\{3,5,6\}\}$.
   Note that all the 4 elements in $\{1,2,3,4\}$ belong to at least
3  sets among the 5 sets in $\{\{1,2,3\},\{1,4,5\},\linebreak  \{2,4,6\},\{3,5,6\}, \{1,2,3,4\}\}$. Then we get that in this case, all the 4 elements in $\{1,2,3,4\}$ belong to at least half of the sets
  in $\cal{G}_3\cup \cal{G}_4$.

 \bigskip
  \quad (c.2) For any  $m\in  \{j_{2l+1},\ldots,j_{n_3-2k}\}, G_m\nsubseteq \{1,2,3,4\}$.   Then $\cal{H}\subset \cal{H}_5\cup \cal{H}_6\cup \cal{H}_{5,6}$,
  where
  \begin{eqnarray*}
  &&\cal{H}_5=\{\{\{1,2,5\},\{1,3,5\}, \{1,4,5\}, \{2,3,5\}, \{2,4,5\},\{3,4,5\} \},\\
  &&\cal{H}_6=\{\{\{1,2,6\},\{1,3,6\}, \{1,4,6\}, \{2,3,6\}, \{2,4,6\},\{3,4,6\} \},\\
  &&\cal{H}_{5,6}=\{\{1,5,6\}, \{2,5,6\}, \{3,5,6\}, \{4,5,6\}\}.
  \end{eqnarray*}
 By (iv), we can easily get that
 $$
 |\cal{H}\cap \cal{H}_5|\leq 2, \  |\cal{H}\cap \cal{H}_6|\leq 2, \  |\cal{H}\cap \cal{H}_{5,6}|\leq 1.
 $$

 \quad  (c.2.1)  $\cal{H}\cap \cal{H}_{5,6}=\emptyset$. Without loss of generality, we assume that $\{1,2,5\}\in \cal{H}$ and $G_{j_{2l+1}}=\{1,2,5\}$. Then by (iv), we know that for any $m\in \{j_{2l+2}, \ldots,j_{n_3-2k}\},\linebreak  G_m\in \{\{3,4,5\}, \{1,3,6\}, \{1,4,6\}, \{2,3,6\}, \{2,4,6\}\}$.
 By (iv) again,  we need only to consider  the following 14 cases:\smallskip

  \quad\quad  (c.2.1.1)  $\cal{H}=\{\{1,2,5\}, \{3,4,5\}\}$.

  \quad\quad (c.2.1.2)  $\cal{H}=\{\{1,2,5\}, \{1,3,6\}\}$.

  \quad\quad  (c.2.1.3)  $\cal{H}=\{\{1,2,5\}, \{1,4,6\}\}$.

 \quad\quad  (c.2.1.4)  $\cal{H}=\{\{1,2,5\}, \{2,3,6\}\}$.

 \quad\quad  (c.2.1.5)  $\cal{H}=\{\{1,2,5\}, \{2,4,6\}\}$.

 \quad\quad  (c.2.1.6)  $\cal{H}=\{\{1,2,5\}, \{3,4,5\}, \{1,3,6\}\}$.

  \quad\quad  (c.2.1.7)  $\cal{H}=\{\{1,2,5\}, \{3,4,5\}, \{1,4,6\}\}$.

   \quad\quad (c.2.1.8)  $\cal{H}=\{\{1,2,5\}, \{3,4,5\}, \{2,3,6\}\}$.

   \quad\quad  (c.2.1.9)  $\cal{H}=\{\{1,2,5\}, \{3,4,5\}, \{2,4,6\}\}$.

  \quad\quad  (c.2.1.10)  $\cal{H}=\{\{1,2,5\}, \{1,3,6\}, \{2,4,6\}\}$.

  \quad\quad  (c.2.1.11)  $\cal{H}=\{\{1,2,5\}, \{1,4,6\}, \{2,3,6\}\}$.

   \quad\quad  (c.2.1.12)  $\cal{H}=\{\{1,2,5\}, \{1,3,6\}, \{2,4,6\}\}$.

   \quad\quad  (c.2.1.13)  $\cal{H}=\{\{1,2,5\}, \{3,4,5\}, \{1,3,6\}, \{2,4,6\}\}$.

   \quad\quad  (c.2.1.14)  $\cal{H}=\{\{1,2,5\}, \{3,4,5\}, \{1,4,6\}, \{2,3,6\}\}$.

\smallskip

\quad  As to  (c.2.1.1), all the 4 elements in $\{1,2,3,4\}$ belong to at least
  two sets among the three sets in $\{\{1,2,5\},\{3,4,5\}, \{1,2,3,4\}\}$. Then we get that all the 4 elements in $\{1,2,3,4\}$ belong to at least half of the sets
  in $\cal{G}_3\cup \cal{G}_4$.
  \smallskip

\quad As to  (c.2.1.2) and (c.2.1.4), all the 3 elements in $\{1,2,3\}$ belong to at least half of the sets
  in $\cal{G}_3\cup \cal{G}_4$.
\smallskip

\quad As to  (c.2.1.3) and (c.2.1.5), all the 3 elements in $\{1,2,4\}$ belong to at least half of the sets
  in $\cal{G}_3\cup \cal{G}_4$.
\smallskip

\quad  As to  (c.2.1.6), all the 4 elements in $\{1,2,3,4\}$ belong to at least
  two sets among the four sets in $\{\{1,2,5\}, \{3,4,5\}, \{1,3,6\},\{1,2,3,4\}\}$. Then we get that  all the 4 elements in $\{1,2,3,4\}$ belong to at least half of the sets
  in $\cal{G}_3\cup \cal{G}_4$.
  \smallskip

\quad As to  (c.2.1.7)-(c.2.1.12), it is easy to check  that  all the 4 elements in $\{1,2,3,4\}$ belong to  at least   two sets among the four sets in $\cal{H}\cup \{\{1,2,3,4\}\}$ and thus
  belong to at least half of the sets  in $\cal{G}_3\cup \cal{G}_4$.
  \smallskip

\quad As to  (c.2.1.13)  and (c.2.1.14), it is easy to check  that  all the 4 elements in $\{1,2,3,4\}$ belong to  at least   three sets among the 5 sets in $\cal{H}\cup \{\{1,2,3,4\}\}$ and thus    belong to at least half of the sets  in $\cal{G}_3\cup \cal{G}_4$.

 \bigskip

 \quad (c.2.2)  $\cal{H}\cap \cal{H}_{5,6}\neq \emptyset$. Without loss of generality, we assume that $\{1,5,6\}\in \cal{H}$ and $G_{j_{2l+1}}=\{1,5,6\}$. Then by (iv), we know that for any $m\in \{j_{2l+2}, \ldots,j_{n_3-2k}\}, G_m\in \{\{2,3,5\}, \{2,4,5\}, \{3,4,5\}, \{2,3,6\}, \{2,4,6\}, \{3,4,6\}\}$.
 By (iv) again,  we need only to consider  the following 3 cases:
 \smallskip

  \quad\quad  (c.2.2.1)  $|\cal{H}\cap \{\{2,3,5\}, \{2,4,5\}, \{3,4,5\}\}|=1,|\cal{H}\cap \{\{2,3,6\}, \{2,4,6\},\linebreak  \{3,4,6\}\}|=0$.

  \quad\quad  (c.2.2.2)  $|\cal{H}\cap \{\{2,3,5\}, \{2,4,5\}, \{3,4,5\}\}|=0,|\cal{H}\cap \{\{2,3,6\}, \{2,4,6\},\linebreak \{3,4,6\}\}|=1$.

  \quad\quad  (c.2.2.3)  $|\cal{H}\cap \{\{2,3,5\}, \{2,4,5\}, \{3,4,5\}\}|=1,|\cal{H}\cap \{\{2,3,6\}, \{2,4,6\},\linebreak \{3,4,6\}\}|=1$.

\smallskip

 \quad As to (c.2.2.1), take $\cal{H}\cap \{\{2,3,5\}, \{2,4,5\}, \{3,4,5\}\}=\{\{2,3,5\}\}$ for example. Now $\cal{H}=\{\{1,5,6\}, \{2,3,5\}\}$. Now all the 3 elements in $\{1,2,3\}$ belong to at least two sets among the  three sets in $\cal{H}\cup \{\{1,2,3,4\}\}$ and thus belong to at least half of the sets in $\cal{G}_3\cup \cal{G}_4$.
 \smallskip

  \quad As to (c.2.2.2), take $\cal{H}\cap \{\{2,3,6\}, \{2,4,6\}, \{3,4,6\}\}=\{\{2,3,6\}\}$ for example. Now $\cal{H}=\{\{1,5,6\}, \{2,3,6\}\}$. Now all the 3 elements in $\{1,2,3\}$ belong to at least two sets among the  three sets in $\cal{H}\cup \{\{1,2,3,4\}\}$ and thus belong to at least half of the sets in $\cal{G}_3\cup \cal{G}_4$.
 \smallskip

  \quad As to (c.2.2.3), take $\cal{H}=\{\{1,2,5\}, \{2,3,5\}, \{2,4,6\}\}$ for example. Now all the 4 elements in $\{1,2,3,4\}$ belong to at least two sets among the  four sets in $\cal{H}\cup \{\{1,2,3,4\}\}$ and thus belong to at least half of the sets in $\cal{G}_3\cup \cal{G}_4$.

\end{itemize}
\end{itemize}

\bigskip

\item[(2.2)]  $n_4\geq 2$.  For any $i,j=1,\ldots,n_4,i\neq j$, we have $H_i\cup H_j=M_6$ or $H_i\cup H_j\in \cal{G}_5$.

\begin{itemize}
\item[(2.2.1)] $n_4$ is an even number and there exists a permutation $(i_1,\ldots,i_{n_4})$ of $(1,\ldots,n_4)$ such that $H_{i_1}\cup H_{i_2}=\cdots=H_{i_{n_4-1}}\cup H_{i_{n_4}}=M_6$. Then all the 6 elements in $M_6$ belong to at least half of the sets in $\cal{G}_4$. Hence it is enough to show that there exist 3 elements in $M_6$ which belong to at least half of the sets in $\cal{G}_3$ or $\cal{G}_3\cup \cal{G}_5$.

\begin{itemize}
\item[(2.2.1.1)] $n_3=1$. Now $\cal{G}_3=\{G_1\}$, and  all  the 3 elements in $G_1$ satisfy the condition.

\item[(2.2.1.2)] $n_3=2$. Now $\cal{G}_3=\{G_1,G_2\}$ and  all  the 3 elements in $G_1\cup G_2$ belong to at least one of the two sets in $\cal{G}_3$.

\item[(2.2.1.3)] $n_3\geq 3$. For any $i,j=1,\ldots,n_3,i\neq j$, we have $G_i\cup G_j=M_6$ or $G_i\cup G_j\in \cal{G}_4\cup \cal{G}_5$.
\smallskip

\quad (a) $n_3$ is an even number and there exists a permutation $(j_1,\ldots,j_{n_3})$ of $(1,\ldots,n_3)$ such that $G_{j_1}\cup G_{j_2}=\cdots=G_{j_{n_3-1}}\cup G_{j_{n_3}}=M_6$.
Then all the 6 elements in $M_6$ belong to half of the sets in $\cal{G}_3$.
\smallskip

\quad (b) $n_3$ is an odd number and there exists a permutation $(j_1,\ldots,j_{n_3})$ of $(1,\ldots,n_3)$ such that $G_{j_1}\cup G_{j_2}=\cdots=G_{j_{n_3-2}}\cup G_{j_{n_3-1}}=M_6$.
Then all the 6 elements in $M_6$ belong to half of the sets in $\{G_{j_1},\ldots,G_{j_{n_3-1}}\}$.
Hence all the 3 elements in $G_{j_{n_3}}$ belong to at least half of the sets in $\cal{G}_3$.
\smallskip

\quad (c) we can decompose $\cal{G}_3$ into two disjoint parts $\{G_{j_1},\ldots,G_{j_{2k}}\}$ and \linebreak $\{G_{j_{2k+1}}, \ldots,G_{j_{n_3}}\}$, where $\{j_1,\ldots,j_{n_3}\}=\{1,\ldots,n_3\}, n_3-2k\geq 2$, and

\quad (i) $G_{j_1}\cup G_{j_2}=\cdots=G_{j_{2k-1}}\cup G_{j_{2k}}=M_6$;

\quad (ii) for any two different indexes $\{i,j\}$ from $\{j_{2k+1},\ldots,j_{n_3}\}$, $G_i\cup G_j\in \cal{G}_4\cup \cal{G}_5$.

Then all the 6 element in $M_6$ belong to half of the sets in $\{G_{j_1},\ldots,G_{j_{2k}}\}$.
Without loss of generality, we assume that $G_{j_{2k+1}}=\{1,2,3\}$. For any $l=j_{2k+2},\ldots,j_{n_3}$, we
 have $G_l\in \{\{1,2,4\}, \{1,2,5\}, \{1,2,6\}, \{1,3,4\}, \{1,3,5\},\linebreak  \{1,3,6\}, \{1,4,5\}, \{1,4,6\}, \{1,5,6\}, \{2,3,4\}, \{2,3,5\}, \{2,3,6\}, \{2,4,5\}, \{2,4,6\}, \linebreak
  \{2,5,6\},  \{3,4,5\}, \{3,4,6\}, \{3,5,6\}\}$. For simplicity, denote $\cal{H}=\{G_{j_{2k+1}}, \ldots,\linebreak G_{j_{n_3}}\}$. By (ii), we know that $|\cal{H}|\leq 10$. Then we have the following 9 cases:
\smallskip

\quad\quad (c.1) $|\cal{H}|=2$. Now $\cal{H}=\{G_{j_{2k+1}},G_{j_{n_3}}\}$ and all the elements in $G_{j_{2k+1}}\cup G_{j_{n_3}}$ belong to at least one of the two sets in $\cal{H}$ and thus belong to at least half of the sets in $\cal{G}_3$.

\quad\quad (c.2) $|\cal{H}|=3$. \quad\quad (c.3) $|\cal{H}|=4$.   \quad\quad (c.4) $|\cal{H}|=5$. \quad\quad (c.5) $|\cal{H}|=6$.

\quad\quad (c.6) $|\cal{H}|=7$. \quad\quad (c.7) $|\cal{H}|=8$. \quad\quad (c.8) $|\cal{H}|=9$.
\quad\quad (c.9) $|\cal{H}|=10$.

\quad As to the cases (c.2)-(c.9), we only give the proof for (c.2), the proofs for the other cases are similar. We omit the details.

\quad As to (c.2),  $\cal{H}=\{G_{j_{2k+1}},G_{j_{2k+2}}, G_{j_{n_3}}\}$ and we have the following 4 cases:
\bigskip

\quad (c.2.1) For any 2-element subset $\{i,j\}$ of $\{j_{2k+1},j_{2k+2}, j_{n_3}\}$, $G_i\cup G_j\in \cal{G}_4$ and $G_{j_{2k+1}}\cup G_{j_{2k+2}}\cup G_{j_{n_3}}\in \cal{G}_4$.
Take $\cal{H}=\{\{1,2,3\}, \{1,2,4\}, \{2,3,4\}\}$ for example. Now by Lemma \ref{lem-1.2}, we know that all the 4 elements in $\{1,2,3,4\}$ belong to at least 2 sets in $\cal{H}$ and thus belong to at least half of the sets in $\cal{G}_3$.
\bigskip

\quad (c.2.2) For any 2-element subset $\{i,j\}$ of $\{j_{2k+1},j_{2k+2}, j_{n_3}\}$, $G_i\cup G_j\in \cal{G}_4$ and $G_{j_{2k+1}}\cup G_{j_{2k+2}}\cup G_{j_{n_3}}\in \cal{G}_5$.
Take $\cal{H}=\{\{1,2,3\}, \{1,2,4\}, \{1,2,5\}\}$ for example. Now $\{1,2,3,4,5\}\in \cal{G}_5$ and thus we have the following 3 cases:
\smallskip

\quad (c.2.2.1) $n_5=1$. Then $\cal{G}_5=\{\{1,2,3,4,5\}\}$. Now all the 5 elements in $\{1,2,3,4,5\}$ belong to at least two sets among the 4 sets in $\cal{H}\cup \cal{G}_5$ and
thus belong to at least half of the sets in $\cal{G}_3\cup \cal{G}_5$. (In fact, by the assumption that $n_5\geq 2$, we don't need consider this case. We write it here for the analysis in the following (c.2.2.3).)
\smallskip

\quad (c.2.2.2) $n_5=2$.  Take $\cal{G}_5=\{\{1,2,3,4,5\}, \{1,2,3,4,6\}\}$ for example. Now it easy to check that all the 4 elements in $\{1,2,3,4\}$ belong to at least 3 sets among the 5 sets in $\cal{H}\cup \cal{G}_5$ and thus belong to at least half of the sets in $\cal{G}_3\cup \cal{G}_5$.
\smallskip

\quad  (c.2.2.3) $n_5\geq 3$. Now by Lemma \ref{lem-1.2}, we know that all the 6 elements in $M_6$ belong to at least $n_5-2$ set(s) in $\cal{G}_5\backslash \{\{1,2,3,4,5\}\}$ and thus belong to at least half of the sets in  $\cal{G}_5\backslash \{\{1,2,3,4,5\}\}$. Then by (c.2.2.1), we know that now all the 5 elements in $\{1,2,3,4,5\}$ belong to at least half of the sets in $\cal{G}_3\cup \cal{G}_5$.
\bigskip

\quad (c.2.3) For any 2-element subset $\{i,j\}$ of $\{j_{2k+1},j_{2k+2}, j_{n_3}\}$, $G_i\cup G_j\in \cal{G}_5$.  Take $\cal{H}=\{\{1,2,3\}, \{1,4,5\}, \{2,5,6\}\}$ for example.
Now $\{\{1,2,3,4,5\}, \linebreak \{1,2,3,5,6\}, \{1,2,4,5,6\}\}\subset \cal{G}_5$ and thus we have the following 3 cases:
\smallskip

\quad (c.2.3.1) $n_5=3$. Then $\cal{G}_5=\{\{1,2,3,4,5\}, \{1,2,3,5,6\}, \{1,2,4,5,6\}\}$. Now all the 6 elements in $M_6$ belong to at least 3 sets among the 6 sets in $\cal{H}\cup \cal{G}_5$ and thus belong to at least half of the sets in $\cal{G}_3\cup \cal{G}_5$.
\smallskip

\quad (c.2.3.2) $n_5=4$. Without loss of generality, we assume that  $\{1,2,3,4,6\}\in \cal{G}_5$. Then by the analysis in (c.2.3.1), we know that all the 5 elements in $\{1,2,3,4,6\}$ belong to at least 4 sets among the 7 sets in $\cal{H}\cup \cal{G}_5$ and thus belong to at least half of the sets in $\cal{G}_3\cup \cal{G}_5$.\smallskip

\quad (c.2.3.2) $n_5\geq 5$.  By Lemma \ref{lem-1.2}, we know that all the 6 elements in $M_6$ belong to at least $|\cal{G}_5\backslash \{\{1,2,3,4,5\}, \{1,2,3,5,6\}, \{1,2,4,5,6\}\}|-1$ set(s) in $\cal{G}_5\backslash \{\{1,2,3,4,5\}, \{1,2,3,5,6\}, \{1,2,4,5,6\}\}$ and thus belong to at least half of the sets in $\cal{G}_5\backslash \{\{1,2,3,4,5\}, \{1,2,3,5,6\}, \{1,2,4,5,6\}\}$. Hence in this case, by (c.2.3.1), we know that all the 6 elements in $M_6$ belong to at least half of the sets in $\cal{G}_3\cup \cal{G}_5$.
\bigskip

\quad (c.2.4) $|\{\{i,j\}|\{i,j\}\subset \{j_{2k+1},j_{2k+2}, j_{n_3}\}, G_i\cup G_j\in \cal{G}_5\}|=2$. Take $\cal{H}=\{\{1,2,3\}, \{1,2,4\}, \{1,5,6\}\}$ for example. Now $\{\{1,2,3,5,6\},\{1,2,4,5,6\}\}\subset \cal{G}_5$. By following the analysis in (c.2.2) and (c.2.3), we can get that there exist at least 3 elements in $M_6$ which belong to at least half of the sets in $\cal{G}_3\cup \cal{G}_5$.
\bigskip

\quad (c.2.5) $|\{\{i,j\}|\{i,j\}\subset \{j_{2k+1},j_{2k+2}, j_{n_3}\}, G_i\cup G_j\in \cal{G}_5\}|=1$. Take $\cal{H}=\{\{1,2,3\}, \{1,2,5\}, \{2,5,6\}\}$ for example. Now $\{1,2,3,5,6\}\in \cal{G}_5$. By following the analysis in (c.2.2) and (c.2.3), we can get that there exist at least 3 elements in $M_6$ which belong to at least half of the sets in $\cal{G}_3\cup \cal{G}_5$.\smallskip
\end{itemize}

\smallskip

\item[(2.2.2)] $n_4$ is an odd number and there exists a permutation $(i_1,\ldots,i_{n_4})$ of $(1,\ldots,n_4)$ such that $H_{i_1}\cup H_{i_2}=\cdots=H_{i_{n_4-2}}\cup H_{i_{n_4-1}}=M_6$. Then all the 6 elements in $M_6$ belong to at least half of the sets in $\{H_{i_1},\ldots,H_{i_{n_4-1}}\}$. Hence it is enough to show that there exist 3 elements in $M_6$
which belong to at least half of the sets in $\{H_{i_{n_4}}\}\cup \cal{G}_3$ or $\{H_{i_{n_4}}\}\cup\cal{G}_3\cup \cal{G}_5$.  Without loss of generality, we assume that $H_{i_{n_4}}=\{1,2,3,4\}$.

\begin{itemize}
\item[(2.2.2.1)] $n_3=1$. Now $\cal{G}_3=\{G_1\}$, and  all  the  elements in $H_{i_{n_4}}\cup G_1$ belong to at least one of the two sets in $\{H_{i_{n_4}}\}\cup \cal{G}_3$.

\item[(2.2.2.2)] $n_3\geq 2$.  For any $i,j=1,\ldots,n_3,i\neq j$, we have $G_i\cup G_j=M_6$ or $G_i\cup G_j\in \cal{G}_4\cup \cal{G}_5$.
\smallskip

\quad (a) $n_3$ is an even number and there exists a permutation $(j_1,\ldots,j_{n_3})$ of $(1,\ldots,n_3)$ such that $G_{j_1}\cup G_{j_2}=\cdots=G_{j_{n_3-1}}\cup G_{j_{n_3}}=M_6$.
Then all the 4 elements in $H_{i_{n_4}}$ belong to at least half of the sets in $\{H_{i_{n_4}}\}\cup \cal{G}_3$.
\smallskip

\quad (b) $n_3$ is an odd number and there exists a permutation $(j_1,\ldots,j_{n_3})$ of $(1,\ldots,n_3)$ such that $G_{j_1}\cup G_{j_2}=\cdots=G_{j_{n_3-2}}\cup G_{j_{n_3-1}}=M_6$.
Then all the 6 elements in $M_6$ belong to half of the sets in $\{G_{j_1},\ldots,G_{j_{n_3-1}}\}$.
Hence all the  elements in $H_{i_{n_4}}\cup G_{j_{n_3}}$ belong to at least half of the sets in $\{H_{i_{n_4}}\}\cup \cal{G}_3$.
\smallskip

\quad (c) we can decompose $\cal{G}_3$ into two disjoint parts $\{G_{j_1},\ldots,G_{j_{2k}}\}$ and \linebreak $\{G_{j_{2k+1}}, \ldots,G_{j_{n_3}}\}$, where $\{j_1,\ldots,j_{n_3}\}=\{1,\ldots,n_3\}, n_3-2k\geq 2$, and

\quad (i) $G_{j_1}\cup G_{j_2}=\cdots=G_{j_{2k-1}}\cup G_{j_{2k}}=M_6$;

\quad (ii) for any two different indexes $\{i,j\}$ from $\{j_{2k+1},\ldots,j_{n_3}\}$, $G_i\cup G_j\in \cal{G}_4\cup \cal{G}_5$.

Then all the 6 element in $M_6$ belong to half of the sets in $\{G_{j_1},\ldots,G_{j_{2k}}\}$. We have the following two cases:
\smallskip

\quad (c.1) There exists $m\in \{j_{2k+1},\ldots,j_{n_3}\}$ such that $G_m\cup \{1,2,3,4\}=M_6$. Without loss of generality, we assume that $G_{j_{2k+1}}=\{1,5,6\}$.
Then for any $l=j_{2k+2},\ldots,j_{n_3}$, we have $G_k\in \{\{1,2,4\}, \{1,2,5\}, \{1,2,6\}, \{1,3,4\}, \{1,3,5\},\linebreak  \{1,3,6\}, \{1,4,5\}, \{1,4,6\},  \{2,3,5\}, \{2,3,6\}, \{2,4,5\}, \{2,4,6\},  \{2,5,6\},  \{3,4,5\},\linebreak \{3,4,6\}, \{3,5,6\}\}$.
\smallskip

\quad (c.2) For any $m\in \{j_{2k+1},\ldots,j_{n_3}\}, G_m\cup \{1,2,3,4\}\neq M_6$. Then
$\{G_{j_{2k+1}}, \ldots, \linebreak G_{j_{n_3}}\}\subset \{\{1,2,3\}, \{1,2,4\}, \{1,2,5\}, \{1,2,6\}, \{1,3,4\}, \{1,3,5\}, \{1,3,6\}, \{1,4,5\},\linebreak \{1,4,6\}, \{2,3,4\}, \{2,3,5\}, \{2,3,6\}, \{2,4,5\}, \{2,4,6\},  \{3,4,5\}, \{3,4,6\}\}$.

\smallskip

\quad As to (c.1) and (c.2), by following the analysis in (2.2.1), we get that there exist at least 3 elements in $M_6$ which belong to at least half of the sets in $\{H_{i_{n_4}}\}\cup\cal{G}_3\cup \cal{G}_5$. We omit the details.

\end{itemize}

\smallskip

\item[(2.2.3)]  We can decompose $\cal{G}_4$ into two disjoint parts $\{H_{i_1},\ldots,H_{i_{2k}}\}$ and $\{H_{i_{2k+1}},\ldots,H_{i_{n_4}}\}$, where $\{i_1,\ldots,i_{n_4}\},n_4-2k\geq 2$, and

\quad (i) $H_{i_1}\cup H_{i_2}=\cdots=H_{i_{n_4-1}}\cup H_{i_{n_4}}=M_6$;

\quad (ii) for any two different indexes $\{i,j\}$ from $\{i_{2k+1},\ldots,i_{n_4}\}$, $H_i\cup H_j\in \cal{G}_5$.

Then all the 6 elements in $M_6$ belong to at least half of the sets in $\{H_{i_1},\ldots,H_{i_{2k}}\}$. Without loss of generality, we assume that $H_{i_{2k+1}}=\{1,2,3,4\}$. Then by (ii), we get that for any $j=i_{2k+2},\ldots,i_{n_4}$,
$H_j\in \cal{H}_5\cup \cal{H}_6$, where
\begin{eqnarray*}
&&\cal{H}_5=\{\{1,2,3,5\}, \{1,2,4,5\}, \{1,3,4,5\}, \{2,3,4,5\}\};\\
&&\cal{H}_6=\{\{1,2,3,6\}, \{1,2,4,6\}, \{1,3,4,6\}, \{2,3,4,6\}\}.
\end{eqnarray*}
For simplicity, denote $\cal{H}=\{H_{i_{2k+1}},\ldots,H_{i_{n_4}}\}$. Then we have the following 3 cases:

\smallskip

\begin{itemize}
\item [(2.2.3.1)]  $\cal{H}\cap \cal{H}_5\neq \emptyset, \cal{H}\cap \cal{H}_6=\emptyset$.  Now,  we have the following 4 cases:

\quad\quad (a) $|\cal{H}\cap \cal{H}_5|=1$.       \quad\quad (b) $|\cal{H}\cap \cal{H}_5|=2$.

\quad\quad (c) $|\cal{H}\cap \cal{H}_5|=3$.       \quad\quad (d) $|\cal{H}\cap \cal{H}_5|=4$.

In the following, we only give the proof for (a). The proofs for other  cases are similar. We omit the details. Without loss of generality, we assume that $\cal{H}\cap \cal{H}_5=\{\{1,2,3,5\}\}$ and thus
$\cal{H}=\{\{1,2,3,4\}, \{1,2,3,5\}\}$.

\bigskip

\quad (a.1) $n_3=1$. Now $\cal{G}_3=\{G_1\}$.  Obviously, all the 3 elements in $\{1,2,3\}$ belong to  at least two sets among the three sets in $\cal{H}\cup \cal{G}_3$ and thus belong to at least half of the sets in
$\cal{F}$.

\bigskip

\quad (a.2) $n_3\geq 2$.  For any $i,j=1,\ldots,n_3,i\neq j$, we have $G_i\cup G_j=M_6$ or $G_i\cup G_j\in \cal{G}_4\cup \cal{G}_5$.

\smallskip

\quad (a.2.1) $n_3$ is an even number and there exists a permutation $(j_1,\ldots,j_{n_3})$ of $(1,\ldots,n_3)$ such that $G_{j_1}\cup G_{j_2}=\cdots=G_{j_{n_3-1}}\cup G_{j_{n_3}}=M_6$.
Then all the 6 elements in $M_6$ belong to half of the sets in $\cal{G}_3$.  Now all the 5 elements in $\{1,2,3,4,5\}$ belong to at least one of the two sets in $\cal{H}$ and thus belong to at least half of the sets in $\cal{F}$.
\bigskip

\quad (a.2.2) $n_3$ is an odd number and there exists a permutation $(j_1,\ldots,j_{n_3})$ of $(1,\ldots,n_3)$ such that $G_{j_1}\cup G_{j_2}=\cdots=G_{j_{n_3-2}}\cup G_{j_{n_3-1}}=M_6$.
Then all the 6 elements in $M_6$ belong to half of the sets in $\{G_{j_1},\ldots,G_{j_{n_3-1}}\}$.  Now, all the 3 elements in $\{1,2,3\}$ belong to  at least two sets among the three
 sets in $\cal{H}\cup \{G_{j_{n_3}}\}$ and thus belong to at least half of the sets in $\cal{F}$.
\bigskip

\quad  (a.2.3) we can decompose $\cal{G}_3$ into two disjoint parts $\{G_{j_1},\ldots,G_{j_{2k}}\}$ and \linebreak $\{G_{j_{2k+1}}, \ldots,G_{j_{n_3}}\}$, where $\{j_1,\ldots,j_{n_3}\}=\{1,\ldots,n_3\}, n_3-2k\geq 2$, and

\quad (i) $G_{j_1}\cup G_{j_2}=\cdots=G_{j_{2k-1}}\cup G_{j_{2k}}=M_6$;

\quad (ii) for any two different indexes $\{i,j\}$ from $\{j_{2k+1},\ldots,j_{n_3}\}$, $G_i\cup G_j\in \cal{G}_4\cup \cal{G}_5$.

Then all the 6 element in $M_6$ belong to half of the sets in $\{G_{j_1},\ldots,G_{j_{2k}}\}$.  Hence in this case, it is enough to show that there exist at least 3 elements in $M_6$ which belong to at least half of the sets in
$\cal{H}\cup \{G_{j_{2k+1}}, \ldots,G_{j_{n_3}}\}$ or $\cal{H}\cup \{G_{j_{2k+1}}, \ldots,G_{j_{n_3}}\}\cup \cal{G}_5$, where $\cal{H}=\{\{1,2,3,4\}, \{1,2,3,5\}\}$.
\bigskip

\quad (a.2.3.1) There exists $m\in \{j_{2k+1},\ldots,j_{n_3}\}$ such that $G_m=\{1,2,3\}=\{1,2,3,4\}\cup \{1,2,3,5\}$. Without loss of generality, we assume that $G_{j_{2k+1}}=\{1,2,3\}$.
\bigskip

\quad (a.2.3.2) There exists $m\in \{j_{2k+1},\ldots,j_{n_3}\}$ such that $G_m\subset \{1,2,3,4\}$ but $G_m\nsubseteq \{1,2,3,5\}$ and for any
$m\in \{j_{2k+1},\ldots,j_{n_3}\},G_m\neq \{1,2,3\}$.  Without loss of generality, we assume that $G_{j_{2k+1}}=\{1,2,4\}$.
\bigskip

\quad (a.2.3.3) There exists $m\in \{j_{2k+1},\ldots,j_{n_3}\}$ such that $G_m\subset \{1,2,3,5\}$ but $G_m\nsubseteq \{1,2,3,4\}$ and for any
$m\in \{j_{2k+1},\ldots,j_{n_3}\},G_m\neq \{1,2,3\}$.  Without loss of generality, we assume that $G_{j_{2k+1}}=\{1,2,5\}$.
\bigskip

\quad (a.2.3.4) For any $m\in \{j_{2k+1},\ldots,j_{n_3}\},G_m\nsubseteq \{1,2,3,4\}$ and $G_m\nsubseteq \{1,2,3,5\}$.  Now
$\cal{H}\subset \{\{1,2,6\}, \{1,3,6\}, \{1,4,5\}, \{1,4,6\}, \{1,5,6\}, \{2,3,6\}, \linebreak \{2,4,5\},  \{2,4,6\}, \{2,5,6\},\{3,4,5\}, \{3,4,6\} \{3,5,6\},\{4,5,6\}\}$.
\bigskip

\quad For the above 4 cases, by following the proof in (2.2.1.3), we can get that there exist at least 3 elements in $M_6$ which belong to at least half of the sets in
$\cal{H}\cup \{G_{j_{2k+1}}, \ldots,G_{j_{n_3}}\}\cup \cal{G}_5$ and thus belong to at least half of the sets in $\cal{F}$.

\bigskip

\item[(2.2.3.2)]    $\cal{H}\cap \cal{H}_6\neq \emptyset, \cal{H}\cap \cal{H}_5=\emptyset$.  The proof is similar to (2.2.3.1). We omit the details.

\bigskip

\item[(2.2.3.3)]   $\cal{H}\cap \cal{H}_5\neq \emptyset, \cal{H}\cap \cal{H}_6\neq \emptyset$.   Without loss of generality, we assume that $\{1,2,3,5\}\in \cal{H}\cap \cal{H}_5$,
 then by (ii) we know that $\cal{H}\cap \cal{H}_6=\{\{1,2,3,6\}\}$ and
$\cal{H}\cap \cal{H}_5=\{\{1,2,3,5\}\}$. Now $\cal{H}=\{\{1,2,3,4\},\{1,2,3,5\},\{1,2,3,6\}\}$.   By following the proof for (2.2.3.1), we get that there exist at least 3 elements in $M_6$
 which belong to at least half of the sets in $\cal{F}$.

\end{itemize}

\end{itemize}
\end{itemize}
\end{itemize}

\vskip 0.5cm
{ \noindent {\bf\large Acknowledgments}

\noindent We thank the anonymous referee for providing helpful comments to
improve and clarify the manuscript. We also thank the support of NNSFC (Grant No. 11771309 and 11871184)
and the Fundamental Research Funds for the Central Universities.

\end{document}